\documentclass[11pt]{amsart}
\usepackage{amsmath,amssymb}
\newtheorem{theorem}{Theorem}[section]
\newtheorem{corollary}[theorem]{Corollary}

\newtheorem{proposition}[theorem]{Proposition}

\begin{document}

\title{Sphere Theorems in Geometry}
\author{Simon Brendle and Richard Schoen}
\address{Department of Mathematics \\
                 Stanford University \\
                 Stanford, CA 94305}
\thanks{The first author was partially supported by a Sloan Foundation Fellowship and by NSF grant DMS-0605223. The second author was partially supported by NSF grant DMS-0604960.}
\begin{abstract}
In this paper, we give a survey of various sphere theorems in geometry. These include the topological sphere theorem of Berger and Klingenberg as well as the differentiable version obtained by the authors. These theorems employ a variety of methods, including geodesic and minimal surface techniques as well as Hamilton's Ricci flow. We also obtain here new results concerning complete manifolds with pinched curvature.
\end{abstract}

\maketitle 

\section{The Topological Sphere Theorem}

The sphere theorem in differential geometry has a long history, dating back to a paper by H.E.~Rauch in 1951. In that paper \cite{Rauch}, Rauch posed the question of whether a compact, simply connected Riemannian manifold $M$ whose sectional curvatures lie in the interval $(1,4]$ is necessarily homeomorphic to the sphere. Around 1960, M.~Berger and W.~Klingenberg gave an affirmative answer to this question:

\begin{theorem}[M.~Berger \cite{Berger1}; W.~Klingenberg \cite{Klingenberg}]
\label{topological.sphere.thm}
Let $M$ be a compact, simply connected Riemannian manifold whose sectional curvatures lie in the interval $(1,4]$. Then $M$ is homeomorphic to $S^n$.
\end{theorem}

More generally, Berger \cite{Berger2} proved that a compact, simply connected Riemannian manifold whose sectional curvatures lie in the interval $[1,4]$ is either homeomorphic to $S^n$ or isometric to a compact symmetric space of rank one. 

K.~Grove and K.~Shiohama proved that the upper bound on the sectional curvature can be replaced by a lower bound on the diameter:

\begin{theorem}[K.~Grove, K.~Shiohama \cite{Grove-Shiohama}] 
\label{diameter.sphere.theorem}
Let $M$ be a compact Riemannian manifold with sectional curvature greater than $1$. If the diameter of $M$ is greater than $\pi/2$, then $M$ is homeomorphic to $S^n$.
\end{theorem}

There is an interesting rigidity statement in the diameter sphere theorem. To describe this result, suppose that $M$ is a compact Riemannian manifold with sectional curvature $K \geq 1$ and diameter $\text{\rm diam}(M) \geq \pi/2$. A theorem of D.~Gromoll and K.~Grove \cite{Gromoll-Grove} asserts that $M$ is either homeomorphic to $S^n$, or locally symmetric, or has the cohomology ring of the $16$-dimensional Cayley plane (see also \cite{Grove-survey}). B.~Wilking \cite{Wilking} proved that, in the latter case, $M$ is isometric to the Cayley plane.

\section{Manifolds with positive isotropic curvature}

M.~Micallef and J.D.~Moore have used harmonic map theory to prove a generalization of Theorem \ref{topological.sphere.thm}. In doing so, they introduced a new curvature condition which they called positive isotropic curvature. A Riemannian manifold $M$ is said to have positive isotropic curvature if 
\[R_{1313} + R_{1414} + R_{2323} + R_{2424} - 2 \, R_{1234} > 0\]
for all points $p \in M$ and all orthonormal four-frames $\{e_1,e_2,e_3,e_4\} \subset T_p M$. We say that $M$ has nonnegative isotropic curvature if 
\[R_{1313} + R_{1414} + R_{2323} + R_{2424} - 2 \, R_{1234} \geq 0\]
for all points $p \in M$ and all orthonormal four-frames $\{e_1,e_2,e_3,e_4\} \subset T_p M$.

We next describe an alternative characterization of positive isotropic curvature, which involves complex notation. To that end, we consider the complexified tangent space $T_p^{\mathbb{C}} M = T_pM \otimes_{\mathbb{R}} \mathbb{C}$. A manifold $M$ has nonnegative isotropic curvature if and only if 
\[R(z,w,\bar{z},\bar{w}) \geq 0\] 
for all points $p \in M$ and all vectors $z,w \in T_p^{\mathbb{C}} M$ satisfying $g(z,z) = g(z,w) = g(w,w) = 0$ (cf. \cite{Micallef-Moore}).

The main theorem of Micallef and Moore is a lower bound for the index of harmonic two-spheres. Recall that the Morse index of a harmonic two-sphere is defined as the number of negative eigenvalues of the second variation operator (counted according to their multiplicities).

\begin{proposition}[M.~Micallef, J.D.~Moore \cite{Micallef-Moore}]
\label{complexified.index.form}
Let $u: S^2 \to M$ be a harmonic map from $S^2$ into a Riemannian manifold $M$. We denote by $E = u^* TM$ the pull-back of the tangent bundle of $M$ under $u$, and by $E^{\mathbb{C}} = E \otimes_{\mathbb{R}} \mathbb{C}$ the complexification of $E$. Moreover, let $I$ be the index form associated with the second variation of energy. Then 
\[I(s,\bar{s}) = 4 \int_{S^2} \big | D_{\frac{\partial}{\partial \bar{z}}} s \big |^2 \, dx \, dy - 4 \int_{S^2} R \Big ( \frac{\partial u}{\partial z},s,\frac{\partial u}{\partial \bar{z}},\bar{s} \Big ) \, dx \, dy.\] 
for all sections $s \in \Gamma(E^{\mathbb{C}})$. Here, $z = x+iy$ denotes the complex coordinate on $S^2$.
\end{proposition}

\textit{Proof of Proposition \ref{complexified.index.form}.} 
Let $I: \Gamma(E^{\mathbb{C}}) \times \Gamma(E^{\mathbb{C}}) \to \mathbb{C}$ denote the complexified index form. Then 
\begin{align*} 
I(s,\bar{s}) 
&= \int_{S^2} \big ( \big | D_{\frac{\partial}{\partial x}} s \big |^2 + \big | D_{\frac{\partial}{\partial y}} s \big |^2 \big ) \, dx \, dy \\ 
&- \int_{S^2} \Big ( R \Big ( \frac{\partial u}{\partial x},s,\frac{\partial u}{\partial x},\bar{s} \Big ) + R \Big ( \frac{\partial u}{\partial y},s,\frac{\partial u}{\partial y},\bar{s} \Big ) \Big ) \, dx \, dy 
\end{align*} 
for all $s \in \Gamma(E^{\mathbb{C}})$. We next define 
\[\frac{\partial u}{\partial z} = \frac{1}{2} \, \Big ( \frac{\partial u}{\partial x} - i \, \frac{\partial u}{\partial y} \Big ) \in \Gamma(E^{\mathbb{C}}), \qquad \frac{\partial u}{\partial \bar{z}} = \frac{1}{2} \, \Big ( \frac{\partial u}{\partial x} + i \, \frac{\partial u}{\partial y} \Big ) \in \Gamma(E^{\mathbb{C}}).\] 
Moreover, for each section $s \in \Gamma(E^{\mathbb{C}})$ we define 
\[D_{\frac{\partial}{\partial z}} s = \frac{1}{2} \, \big ( D_{\frac{\partial}{\partial x}} s - i \, D_{\frac{\partial}{\partial y}} s \big ), \qquad D_{\frac{\partial}{\partial \bar{z}}} s = \frac{1}{2} \, \big ( D_{\frac{\partial}{\partial x}} s + i \, D_{\frac{\partial}{\partial y}} s \big ).\] 
With this understood, the complexified index form can be written in the form 
\begin{align*} 
I(s,\bar{s}) 
&= 2 \int_{S^2} \big ( \big | D_{\frac{\partial}{\partial z}} s \big |^2 + \big | D_{\frac{\partial}{\partial \bar{z}}} s \big |^2 \big ) \, dx \, dy \\ 
&- 2 \int_{S^2} \Big ( R \Big ( \frac{\partial u}{\partial \bar{z}},s,\frac{\partial u}{\partial z},\bar{s} \Big ) + R \Big ( \frac{\partial u}{\partial z},s,\frac{\partial u}{\partial \bar{z}},\bar{s} \Big ) \Big ) \, dx \, dy 
\end{align*} 
for all $s \in \Gamma(E^{\mathbb{C}})$. Integration by parts yields 
\begin{align*} 
&\int_{S^2} \big ( \big | D_{\frac{\partial}{\partial z}} s \big |^2 - \big | D_{\frac{\partial}{\partial \bar{z}}} s \big |^2 \big ) \, dx \, dy \\ 
&= \int_{S^2} g \big ( D_{\frac{\partial}{\partial z}} D_{\frac{\partial}{\partial \bar{z}}} s - D_{\frac{\partial}{\partial \bar{z}}} D_{\frac{\partial}{\partial z}} s,\bar{s} \big ) \, dx \, dy \\ 
&= -\int_{S^2} R \Big ( \frac{\partial u}{\partial z},\frac{\partial u}{\partial \bar{z}},s,\bar{s} \Big ) \, dx \, dy \\ 
&= \int_{S^2} \Big ( R \Big ( \frac{\partial u}{\partial \bar{z}},s,\frac{\partial u}{\partial z},\bar{s} \Big ) - R \Big ( \frac{\partial u}{\partial z},s,\frac{\partial u}{\partial \bar{z}},\bar{s} \Big ) \Big ) \, dx \, dy 
\end{align*} 
for all $s \in \Gamma(E^{\mathbb{C}})$. Putting these facts together, the assertion follows. \\

\begin{theorem}[M.~Micallef, J.D.~Moore \cite{Micallef-Moore}]
\label{Morse.index}
Let $u: S^2 \to M$ be a harmonic map from $S^2$ into a Riemannian manifold $M$. If $M$ has positive isotropic curvature, then $u$ has Morse index at least $[\frac{n-2}{2}]$.
\end{theorem}

\textit{Proof of Theorem \ref{Morse.index}.} 
We denote by $E = u^* TM$ the pull-back of the tangent bundle of $M$, and by $E^{\mathbb{C}}$ the complexification of $E$. 
Let $z = x+iy$ the complex coordinate on $S^2$. As above, we define 
\[D_{\frac{\partial}{\partial \bar{z}}} s = \frac{1}{2} \, (D_{\frac{\partial}{\partial x}} s + i \, D_{\frac{\partial}{\partial y}} s)\] 
for each section $s \in \Gamma(E^{\mathbb{C}})$. We say that $s \in \Gamma(E^{\mathbb{C}})$ is holomorphic if $D_{\frac{\partial}{\partial \bar{z}}} s = 0$.

Let $\mathcal{H}$ denote the space of holomorphic sections of $E^{\mathbb{C}}$. Given two holomorphic sections $s_1,s_2 \in \mathcal{H}$, the inner product $g(s_1,s_2)$ defines a holomorphic function on $S^2$. Consequently, the function $g(s_1,s_2)$ is constant. This defines a symmetric bilinear form 
\[\mathcal{H} \times \mathcal{H} \to \mathbb{C}, \qquad (s_1,s_2) \mapsto g(s_1,s_2).\] 

By assumption, the map $u: S^2 \to M$ is harmonic. Hence, $\frac{\partial u}{\partial z}$ is a holomorphic 
section of $E^{\mathbb{C}}$. Since $u$ is smooth at the north pole on $S^2$, the section $\frac{\partial u}{\partial z}$ 
vanishes at the north pole. Thus, we conclude that $g(\frac{\partial u}{\partial z},s) = 0$ for every holomorphic section $s \in \mathcal{H}$. In particular, we have $g(\frac{\partial u}{\partial z},\frac{\partial u}{\partial z}) = 0$.

By the Grothendieck splitting theorem (cf. \cite{Grothendieck}), the bundle $E^{\mathbb{C}}$ splits as a direct sum of holomorphic line subbundles; that is,  
\[E^{\mathbb{C}} = L_1 \oplus L_2 \oplus \hdots \oplus L_n.\] 
We assume that the line bundles $L_1,L_2, \hdots, L_n$ are chosen so that 
\[c_1(L_1) \geq c_1(L_2) \geq \hdots \geq c_1(L_n).\] 
Note that $c_1(L_1), c_1(L_2), \hdots, c_1(L_n)$ are uniquely determined, but $L_1,L_2, \hdots, L_n$ are not. By definition, $E^{\mathbb{C}}$ is the complexification of a real bundle. In particular, the bundle $E^{\mathbb{C}}$ is canonically isomorphic to its dual bundle. From this, we deduce that 
\[c_1(L_k) + c_1(L_{n-k+1}) = 0\] 
for $k = 1, \hdots, n$ (see \cite{Micallef-Moore}, p.~209). 

For each $k \in \{1, \hdots, n\}$, we denote by $F^{(k)}$ the direct sum of all line bundles $L_j$ except $L_k$ and $L_{n-k+1}$. More precisely, we define 
\[F^{(k)} = \bigoplus_{j \in \mathcal{J}^{(k)}} L_j,\] 
where $\mathcal{J}^{(k)} = \{1, \hdots, n\} \setminus \{k,n-k+1\}$. Note that $\bigcap_{k=1}^n F^{(k)} = \{0\}$. Moreover, we have $c_1(F^{(k)}) = 0$ and $\text{\rm rank} \, F^{(k)} \geq n-2$. Let $\mathcal{H}^{(k)} \subset \mathcal{H}$ denote the space of holomorphic sections of $F^{(k)}$. It follows from the Riemann-Roch theorem that $\dim_{\mathbb{C}} \mathcal{H}^{(k)} \geq n-2$.

Fix an integer $k \in \{1, \hdots, n\}$ such that $\frac{\partial u}{\partial z} \notin \Gamma(F^{(k)})$. Since $\dim_{\mathbb{C}} \mathcal{H}^{(k)} \geq n-2$, there exists a subspace $\tilde{\mathcal{H}} \subset \mathcal{H}^{(k)}$ such that $\dim_{\mathbb{C}} \tilde{\mathcal{H}} \geq [\frac{n-2}{2}]$ and $g(s,s) = 0$ for all $s \in \tilde{\mathcal{H}}$. We claim that the resctriction of $I$ to $\tilde{\mathcal{H}}$ is negative definite. To see this, consider a section $s \in \tilde{\mathcal{H}}$. Since $s$ is holomorphic, we have 
\[I(s,\bar{s}) = -4 \int_{S^2} R \Big ( \frac{\partial u}{\partial z},s,\frac{\partial u}{\partial \bar{z}},\bar{s} \Big ) \, dx \, dy\] 
by Proposition \ref{complexified.index.form}. Moreover, we have $g(s,s) = g(\frac{\partial u}{\partial z},s) = g(\frac{\partial u}{\partial z},\frac{\partial u}{\partial z}) = 0$. Since $M$ has positive isotropic curvature, it follows that 
\[R \Big ( \frac{\partial u}{\partial z},s,\frac{\partial u}{\partial \bar{z}},\bar{s} \Big ) \geq 0.\] 
Putting these facts together, we conclude that $I(s,\bar{s}) \leq 0$. It remains to analyze the case of equality. If $I(s,\bar{s}) = 0$, then $s = f \, \frac{\partial u}{\partial z}$ for some meromorphic function $f: S^2 \to \mathbb{C}$. 
However, $\frac{\partial u}{\partial z} \notin \Gamma(F^{(k)})$ by our choice of $k$. Since $s \in \Gamma(F^{(k)})$, it follows that $f$ vanishes identically. Therefore, the restriction of $I$ to $\tilde{\mathcal{H}}$ is negative definite. 

We now complete the proof of Theorem \ref{Morse.index}. Suppose that $m < [\frac{n-2}{2}]$, where $m$ denotes the number of negative eigenvalues of the second variation operator. Then $\dim_{\mathbb{C}} \tilde{\mathcal{H}} > m$. Consequently, there exists a non-vanishing section $s \in \tilde{\mathcal{H}}$ which is orthogonal to the first $m$ eigenfunctions of the second variation operator. Since $s \in \tilde{\mathcal{H}}$, we have $I(s,\bar{s}) < 0$. On the other hand, we have $I(s,\bar{s}) \geq 0$ since $s$ is orthogonal to the first $m$ eigenfunctions of the second variation operator. This is a contradiction. \\

Combining their index estimate with the existence theory of Sacks and Uhlenbeck \cite{Sacks-Uhlenbeck}, Micallef and Moore obtained the following result:

\begin{theorem}[M.~Micallef, J.D.~Moore \cite{Micallef-Moore}]
\label{homotopy.sphere}
Let $M$ be a compact simply connected Riemannian manifold with positive isotropic curvature. Then $M$ is a homotopy sphere. Hence, if $n \geq 4$, then $M$ is homeomorphic to $S^n$.
\end{theorem}

\textit{Sketch of the proof of Theorem \ref{homotopy.sphere}.} 
Suppose that $\pi_j(M) \neq 0$ for some integer $j \geq 2$. By a theorem of Sacks and Uhlenbeck \cite{Sacks-Uhlenbeck}, there exists a harmonic map $u: S^2 \to M$ with Morse index less than $j - 1$. On the other hand, any harmonic map $u: S^2 \to M$ has Morse index at least $[\frac{n}{2}] - 1$ by Theorem \ref{Morse.index}. Putting these facts together, we obtain $j > [\frac{n}{2}]$. Thus, $\pi_j(M) = 0$ for $j = 2, \hdots, [\frac{n}{2}]$. Since $M$ is simply connected, the Hurewicz theorem implies that $\pi_j(M) = 0$ for $j = 1, \hdots, n-1$. Consequently, $M$ is a homotopy sphere. \\

We say that $M$ has pointwise $1/4$-pinched sectional curvatures if $0 < K(\pi_1) < 4 \, K(\pi_2)$ for all points $p \in M$ and all two-planes $\pi_1,\pi_2 \subset T_p M$. It follows from Berger's inequality (see e.g. \cite{Karcher}) that every manifold with pointwise $1/4$-pinched sectional curvatures has positive isotropic curvature. Hence, Theorem \ref{homotopy.sphere} generalizes the classical sphere theorem of Berger and Klingenberg.

The topology of non-simply connected manifolds with positive isotropic curvature is not fully understood. It has been conjectured that the fundamental group of a compact manifold $M$ with positive isotropic curvature is virtually free in the sense that it contains a free subgroup of finite index (see \cite{Fraser2},\cite{Gromov}). A.~Fraser has obtained an important result in this direction:

\begin{theorem}[A.~Fraser \cite{Fraser2}] 
\label{Fraser.thm}
Let $M$ be a compact Riemannian manifold of dimension $n \geq 5$ with positive isotropic curvature. Then the fundamental group of $M$ does not contain a subgroup isomorphic to $\mathbb{Z} \oplus \mathbb{Z}$.
\end{theorem}

The proof of Theorem \ref{Fraser.thm} relies on the existence theory of Schoen and Yau \cite{Schoen-Yau}, and a careful study of the second variation of area (see also \cite{Fraser1},\cite{Schoen-survey}). The proof also uses the following result due to A.~Fraser (see \cite{Fraser2}, Section 3): 

\begin{proposition}[A.~Fraser \cite{Fraser2}]
\label{degree.one.map}
Let $h$ be a Riemannian metric on $T^2$ with the property that every non-contractible loop in $(T^2,h)$ has length at least $1$. Moreover, let $S^2$ be the two-sphere equipped with its standard metric of constant curvature $1$. Then there exists a degree-one map $f$ from $(T^2,h)$ to $S^2$ such that $|Df| \leq C$, where $C$ is a numerical constant.
\end{proposition}

\textit{Proof of Proposition \ref{degree.one.map}.} 
Let $\Sigma$ be the universal cover of $(T^2,h)$, and let $\pi: \Sigma \to (T^2,h)$ denote the covering projection. Note that $\Sigma$ is diffeomorphic to $\mathbb{R}^2$. For each positive integer $k$, there exists a unit-speed geodesic $\gamma_k: [-k,k] \to \Sigma$ such that $d(\gamma_k(k),\gamma_k(-k)) = 2k$. Passing to the limit as $k \to \infty$, we obtain a unit-speed geodesic $\gamma: \mathbb{R} \to \Sigma$ such that $d(\gamma(t_1),\gamma(t_2)) = |t_1 - t_2|$ for all $t_1,t_2 \in \mathbb{R}$. In particular, $\gamma(t_1) \neq \gamma(t_2)$ whenever $t_1 \neq t_2$. By the Jordan curve theorem, the complement $\Sigma \setminus \{\gamma(t): t \in \mathbb{R}\}$ has exactly two connected components, which we denote by $\Omega_1$ and $\Omega_2$. 

We next define functions $D_1: \Sigma \to \mathbb{R}$ and $D_2: \Sigma \to \mathbb{R}$ by 
\[D_1(p) = \begin{cases} \inf \{d(\gamma(t),p): t \in \mathbb{R}\} & \text{\rm for $p \in \Omega_1$} \\ -\inf \{d(\gamma(t),p): t \in \mathbb{R}\} & \text{\rm for $p \in \Omega_2$} \\ 0 & \text{\rm otherwise} \end{cases}\] 
and 
\[D_2(p) = d(\gamma(0),p) - 1.\] 
Clearly, $|D_j(p) - D_j(q)| \leq d(p,q)$ for all points $p,q \in \Sigma$. Let 
\[\mathcal{Q} = \Big \{ p \in \Sigma: D_1(p)^2 + D_2(p)^2 \leq \frac{1}{64} \Big \}.\] 
We claim that 
\[\mathcal{Q} \subset B_{1/3}(\gamma(1)) \cup B_{1/3}(\gamma(-1)).\] 
To see this, we consider a point $p \in \mathcal{Q}$. Then there exists a real number $t$ such 
that $d(\gamma(t),p) = |D_1(p)|$. This implies 
\begin{align*} 
\big | |t| - 1 \big | 
&= \big | d(\gamma(0),\gamma(t)) - 1 \big | \\ 
&\leq \big | d(\gamma(0),\gamma(t)) - d(\gamma(0),p) \big | + \big | d(\gamma(0),p) - 1 \big | \\ 
&\leq d(\gamma(t),p) + \big | d(\gamma(0),p) - 1 \big | \\ 
&= |D_1(p)| + |D_2(p)|. 
\end{align*} From this, we deduce that 
\begin{align*} 
&\min \big \{ d(\gamma(1),p),d(\gamma(-1),p) \big \} \\ 
&\leq \min \big \{ d(\gamma(1),\gamma(t)),d(\gamma(-1),\gamma(t)) \big \} + d(\gamma(t),p) \\ 
&= \big | |t| - 1 \big | + d(\gamma(t),p) \\ 
&\leq 2 \, |D_1(p)| + |D_2(p)| \\ 
&< \frac{1}{3}. 
\end{align*}
Thus, $\mathcal{Q} \subset B_{1/3}(\gamma(1)) \cup B_{1/3}(\gamma(-1))$.

We next define $\mathcal{R} = \mathcal{Q} \cap B_1(\gamma(1))$. Clearly, $\mathcal{R} \subset B_{1/3}(\gamma(1))$. This implies $D_1(p)^2 + D_2(p)^2 = \frac{1}{64}$ for all points $p \in \partial \mathcal{R}$. Hence, the map $(D_1,D_2): \mathcal{R} \to B_{1/8}(0)$ maps $\partial \mathcal{R}$ into $\partial B_{1/8}(0)$. The map $(D_1,D_2)$ is smooth in a neighborhood of $\gamma(1)$. Moreover, the differential of $(D_1,D_2)$ at the point $\gamma(1)$ is non-singular. Since $(D_1(p),D_2(p)) \neq (D_1(\gamma(1)),D_2(\gamma(1)))$ for all $p \in \mathcal{R} \setminus \{\gamma(1)\}$, we conclude that the map $(D_1,D_2): \mathcal{R} \to B_{1/8}(0)$ has degree one.

In the next step, we approximate the functions $D_1$ and $D_2$ by smooth functions. Let $\delta$ be an arbitrary positive real number. Using the convolution procedure of Greene and Wu (see \cite{Greene-Wu1},\cite{Greene-Wu2}), 
we can construct smooth functions $\tilde{D}_1: \mathcal{R} \to \mathbb{R}$ and $\tilde{D}_2: \mathcal{R} \to \mathbb{R}$ such that 
\[|\tilde{D}_j(p) - D_j(p)| \leq \delta\] 
and 
\[|\tilde{D}_j(p) - \tilde{D}_j(q)| \leq 2 \, d(p,q)\] 
for all points $p,q \in \mathcal{R}$. 

Fix a cut-off function $\eta: [0,\infty) \to [0,1]$ such that $\eta(s) = 2$ for $s \leq 2$ and $\eta(s) = 0$ for $s \geq 3$. We define smooth maps $\varphi: \mathbb{R}^2 \to \mathbb{R}^3 \setminus \{0\}$ and $\psi: \mathbb{R}^2 \to S^2$ by 
\[\varphi(x_1,x_2) = \big ( x_1 \, \eta(x_1^2+x_2^2), x_2 \, \eta(x_1^2+x_2^2), 1-x_1^2-x_2^2 \big )\] 
and 
\[\psi(x_1,x_2) = \frac{\varphi(x_1,x_2)}{|\varphi(x_1,x_2)|}.\] 
In particular, $\psi(x_1,x_2) = (0,0,-1)$ whenever $x_1^2 + x_2^2 \geq 3$. We now define a map $F: \mathcal{R} \to S^2$ by 
\[F(p) = \psi \big ( 16 \, \tilde{D}_1(p),16 \, \tilde{D}_2(p) \big ).\] 
There exists a numerical constant $C$ such that $d(F(p),F(q)) \leq C \, d(p,q)$ for all points $p,q \in \mathcal{R}$. 
Moreover, $F$ maps a neighborhood of the boundary $\partial \mathcal{R}$ to the south pole on $S^2$. It is easy to 
see that the map $F: \mathcal{R} \to S^2$ has degree one. 

By assumption, every non-contractible loop in $(T^2,h)$ has length at least $1$. Hence, if $p,q$ are two distinct points in $\Sigma$ satisfying $\pi(p) = \pi(q)$, then $d(p,q) \geq 1$. Since $\mathcal{R} \subset B_{1/3}(\gamma(1))$, it follows that the restriction $\pi |_{\mathcal{R}}$ is injective. We now define a map $f: (T^2,h) \to S^2$ by 
\[f(y) = \begin{cases} F(p) & \text{\rm if $y = \pi(p)$ for some point $p \in \mathcal{R}$} \\ (0,0,-1) & \text{\rm otherwise} \end{cases}\] 
for $y \in T^2$. It is straightforward to verify that $f$ has all the required properties. This completes the proof of Proposition \ref{degree.one.map}. \\

We note that the minimal surface arguments in \cite{Fraser2} can be extended to the case $n = 4$ provided that $M$ is orientable. 

\begin{theorem} 
\label{Fraser.thm.dim.4}
Let $M$ be a compact orientable four-manifold with positive isotropic curvature. Then the fundamental group of $M$ does not contain a subgroup isomorphic to $\mathbb{Z} \oplus \mathbb{Z}$.
\end{theorem}

\textit{Proof of Theorem \ref{Fraser.thm.dim.4}.} 
Suppose that $\pi_1(M)$ contains a subgroup $G$ which is isomorphic to $\mathbb{Z} \oplus \mathbb{Z}$. For each positive integer $k$, we denote by $G_k$ the subgroup of $\pi_1(M)$ corresponding to $k\mathbb{Z} \oplus k\mathbb{Z}$. Moreover, let 
\[\Lambda_k = \inf \{L(\alpha): \text{\rm $\alpha$ is a non-contractible loop in $M$ with $[\alpha] \in G_k$}\}.\]
Note that $\Lambda_k \to \infty$ as $k \to \infty$.

Fix $k$ sufficiently large. By a theorem of Schoen and Yau \cite{Schoen-Yau}, there exists a branched conformal minimal immersion $u: T^2 \to M$ with the property that $u_*: \pi_1(T^2) \to \pi_1(M)$ is injective and maps $\pi_1(T^2)$ to $G_k$. Moreover, the map $u$ minimizes area in its homotopy class. Hence, $u$ is stable. We next consider the normal bundle of the surface $u(T^2)$. We denote by $E$ the pull-back, under $u$, of the normal bundle of $u(T^2)$. Note that $E$ is a smooth vector bundle of rank $2$, even across branch points. (This follows from the analysis of branch points in \cite{Gulliver},\cite{Micallef-White}.) Since $M$ and $T^2$ are orientable, we conclude that $E$ is orientable. Let $E^{\mathbb{C}} = E \otimes_{\mathbb{R}} \mathbb{C}$ be the complexification of $E$. Since $E$ is orientable, the complexified bundle $E^{\mathbb{C}}$ splits as a direct sum of two holomorphic line bundles $E^{(1,0)}$ and $E^{(0,1)}$. Here, $E^{(1,0)}$ consists of all vectors of the form $a (v - iw) \in E^{\mathbb{C}}$, where $a \in \mathbb{C}$ and $\{v,w\}$ is a positively oriented orthonormal basis of $E$. Similarly, $E^{(0,1)}$ consists of all vectors of the form $a (v + iw) \in E^{\mathbb{C}}$, where $a \in \mathbb{C}$ and $\{v,w\}$ is a positively oriented orthonormal basis of $E$. Since $E^{\mathbb{C}}$ is the complexification of a real bundle, we have $c_1(E^{(1,0)}) + c_1(E^{(0,1)}) = c_1(E^{\mathbb{C}}) = 0$. Without loss of generality, we may assume that $c_1(E^{(1,0)}) \geq 0$. (Otherwise, we choose the opposite orientation on $E$.)

Since $u$ is stable, we have 
\begin{equation} 
\label{stability.1}
\int_{T^2} \big | D_{\frac{\partial}{\partial \bar{z}}}^\perp s \big |^2 \, dx \, dy \geq \int_{T^2} R \Big ( \frac{\partial u}{\partial z},s,\frac{\partial u}{\partial \bar{z}},\bar{s} \Big ) \, dx \, dy 
\end{equation}
for all sections $s \in \Gamma(E^{\mathbb{C}})$ (see \cite{Fraser1},\cite{Schoen-survey}). Every section $s \in \Gamma(E^{(1,0)})$ is isotropic, i.e. $g(s,s) = 0$. Since $M$ has positive isotropic curvature, there exists a positive constant $\kappa$ such that  
\[R \Big ( \frac{\partial u}{\partial z},s,\frac{\partial u}{\partial \bar{z}},\bar{s} \Big ) \geq \kappa \, \Big | \frac{\partial u}{\partial z} \Big |^2 \, |s|^2\] 
for all sections $s \in \Gamma(E^{(1,0)})$. Putting these facts together, we obtain 
\begin{equation} 
\label{stability.2} 
\int_{T^2} \big | D_{\frac{\partial}{\partial \bar{z}}}^\perp s \big |^2 \, dx \, dy \geq \kappa \int_{T^2} \Big | \frac{\partial u}{\partial z} \Big |^2 \, |s|^2 \, dx \, dy 
\end{equation}
for all $s \in \Gamma(E^{(1,0)})$. Moreover, we can find a positive constant $\varepsilon = \varepsilon(k)$ such that 
\begin{equation} 
\label{stability.3}
\int_{T^2} \big | D_{\frac{\partial}{\partial \bar{z}}}^\perp s \big |^2 \, dx \, dy + \frac{1}{2} \, \kappa \int_{T^2} \Big | \frac{\partial u}{\partial z} \Big |^2 \, |s|^2 \, dx \, dy \geq \frac{1}{2} \, \kappa \, \varepsilon \int_{T^2} |s|^2 \, dx \, dy 
\end{equation} 
for all $s \in \Gamma(E^{(1,0)})$. Taking the arithmetic mean of (\ref{stability.2}) and (\ref{stability.3}), we obtain 
\begin{equation} 
\label{stability.4}
\int_{T^2} \big | D_{\frac{\partial}{\partial \bar{z}}}^\perp s \big |^2 \, dx \, dy 
\geq \frac{1}{4} \, \kappa \int_{T^2} \Big ( \Big | \frac{\partial u}{\partial z} \Big |^2 + \varepsilon \Big ) \, |s|^2 \, dx \, dy 
\end{equation}
for all $s \in \Gamma(E^{(1,0)})$. 

We next define a Riemannian metric $h$ on $T^2$ by 
\[h = u^* g + 2\varepsilon \, (dx \otimes dx + dy \otimes dy) = u^* g + \varepsilon \, (dz \otimes d\bar{z} + d\bar{z} \otimes dz).\] 
Every non-contractible loop in $(T^2,h)$ has length at least $\Lambda_k$. By Proposition \ref{degree.one.map}, there exists a degree-one map $f$ from $(T^2,h)$ to the standard sphere $S^2$ such that $\Lambda_k \, |Df| \leq C$. This implies 
\begin{equation} 
\label{bound.for.df}
\Lambda_k^2 \, \Big | \frac{\partial f}{\partial z} \Big |^2 \leq C_1 \, \Big | \frac{\partial}{\partial z} \Big |_h^2 = C_1 \, \Big ( \Big | \frac{\partial u}{\partial z} \Big |^2 + \varepsilon \Big ),
\end{equation} 
where $C_1$ is a positive constant independent of $k$.

Fix a holomorphic line bundle $L$ over $S^2$ with $c_1(L) > 0$. We also fix a metric and a connection on $L$. Finally, we fix sections $\omega_1,\omega_2 \in \Gamma(L^*)$ such that $|\omega_1| + |\omega_2| \geq 1$ at each point on $S^2$. 

Let $\xi = f^* L$ be the pull-back of $L$ under the map $f$. Since $f$ has degree one, we have $c_1(\xi) > 0$. Since $c_1(E^{(1,0)}) \geq 0$, it follows that $c_1(E^{(1,0)} \otimes \xi) > 0$. By the Riemann-Roch theorem, the bundle $E^{(1,0)} \otimes \xi$ admits a non-vanishing holomorphic section, which we denote by $\sigma$. For $j = 1,2$, we define $\tau_j = f^*(\omega_j) \in \Gamma(\xi^*)$ and $s_j = \sigma \otimes \tau_j \in \Gamma(E^{(1,0)})$. Since $\sigma$ is holomorphic, we have 
\[D_{\frac{\partial}{\partial \bar{z}}}^\perp s_j = \sigma \otimes \nabla_{\frac{\partial}{\partial \bar{z}}} \tau_j,\] 
where $\nabla$ denotes the connection on $\xi^*$. We next observe that 
\[\big | \nabla_{\frac{\partial}{\partial \bar{z}}} \tau_j \big |^2 = \big | \nabla_{\frac{\partial f}{\partial \bar{z}}} \omega_j \big |^2 \leq C_2 \, \Big | \frac{\partial f}{\partial z} \Big |^2,\] 
where $C_2$ is a positive constant independent of $k$. This implies 
\[\big | D_{\frac{\partial}{\partial \bar{z}}}^\perp s_j \big |^2 = \big | \nabla_{\frac{\partial}{\partial \bar{z}}} \tau_j \big |^2 \, |\sigma|^2 \leq C_2 \, \Big | \frac{\partial f}{\partial z} \Big |^2 \, |\sigma|^2\] 
for $j = 1,2$. Using (\ref{bound.for.df}), we obtain 
\[\Lambda_k^2 \, \big | D_{\frac{\partial}{\partial \bar{z}}}^\perp s_j \big |^2 \leq C_1C_2 \, \Big ( \Big | \frac{\partial u}{\partial z} \Big |^2 + \varepsilon \Big ) \, |\sigma|\] 
for $j = 1,2$. From this, we deduce that 
\begin{align} 
&\Lambda_k^2 \int_{T^2} \big ( \big | D_{\frac{\partial}{\partial \bar{z}}}^\perp s_1 \big |^2 + \big | D_{\frac{\partial}{\partial \bar{z}}}^\perp s_2 \big |^2 \big ) \, dx \, dy \notag \\ 
&\leq 2 \, C_1C_2 \int_{T^2} \Big ( \Big | \frac{\partial u}{\partial z} \Big |^2 + \varepsilon \Big ) \, |\sigma|^2 \, dx \, dy. 
\end{align} 
Note that 
\[|s_1| + |s_2| = |\sigma| \, (|\tau_1| + |\tau_2|) \geq |\sigma|\] 
at each point on $T^2$. Hence, it follows from (\ref{stability.4}) that 
\begin{align} 
&\int_{T^2} \big ( \big | D_{\frac{\partial}{\partial \bar{z}}}^\perp s_1 \big |^2 + \big | D_{\frac{\partial}{\partial \bar{z}}}^\perp s_1 \big |^2 \big ) \, dx \, dy \notag \\ 
&\geq \frac{1}{4} \, \kappa \int_{T^2} \Big ( \Big | \frac{\partial u}{\partial z} \Big |^2 + \varepsilon \Big ) \, (|s_1|^2 + |s_2|^2) \, dx \, dy \\ 
&\geq \frac{1}{8} \, \kappa \int_{T^2} \Big ( \Big | \frac{\partial u}{\partial z} \Big |^2 + \varepsilon \Big ) \, |\sigma|^2 \, dx \, dy. \notag 
\end{align} 
Thus, we conclude that $\kappa \, \Lambda_k^2 \leq 16 \, C_1C_2$. This contradicts the fact that $\Lambda_k \to \infty$ as $k \to \infty$. \\

In the remainder of this section, we describe sufficient conditions for the vanishing of the second Betti number. M.~Berger \cite{Berger2} proved that the second Betti number of a manifold with pointwise $1/4$-pinched sectional curvatures is equal to $0$. In even dimensions, the same result holds under the weaker assumption that $M$ has positive isotropic curvature: 

\begin{theorem}[M.~Micallef, M.~Wang \cite{Micallef-Wang}] 
\label{Betti.number.even}
Let $M$ be a compact Riemannian manifold of dimension $n \geq 4$. Suppose that $n$ is even and $M$ has positive isotropic curvature. Then the second Betti number of $M$ vanishes.
\end{theorem}

\textit{Proof of Theorem \ref{Betti.number.even}.} 
Suppose that $\psi$ is a non-vanishing harmonic two-form on $M$. It follows from the Bochner formula that 
\[\Delta \psi_{ik} = \sum_{j=1}^n \text{\rm Ric}_i^j \, \psi_{jk} + \sum_{j=1}^n \text{\rm Ric}_k^j \, \psi_{ij} - 2 \sum_{j,l=1}^n R_{ijkl} \, \psi^{jl},\] 
where $\Delta \psi = \sum_{j,l=1}^n g^{jl} \, D_{j,l}^2 \psi$ denotes the rough Laplacian of $\psi$. Fix a point $p \in M$ where the function $|\psi|^2$ attains its maximum. At the point $p$, we have $|\psi|^2 > 0$ and $\Delta(|\psi|^2) \leq 0$. This implies 
\begin{align} 
\label{maximum.point.1}
0 &\geq \Delta \bigg ( \sum_{i,k=1}^n \psi_{ik} \, \psi^{ik} \bigg ) \notag \\ 
&\geq 2 \sum_{i,k=1}^n \Delta \psi_{ik} \, \psi^{ik} = 4 \sum_{i,j,k,l=1}^n (\text{\rm Ric}^{ij} \, g^{kl} - R^{ijkl}) \, \psi_{ik} \, \psi_{jl}
\end{align}
at the point $p$. In order to analyze the curvature term on the right hand side, we write $n = 2m$. We can find an orthonormal basis $\{v_1,w_1,v_2,w_2,\hdots,v_m,w_m\}$ of $T_p M$ and real numbers $\lambda_1,\hdots,\lambda_m$ such that 
\begin{align*} 
&\psi(v_\alpha,w_\beta) = \lambda_\alpha \, \delta_{\alpha\beta} \\ 
&\psi(v_\alpha,v_\beta) = \psi(w_\alpha,w_\beta) = 0 
\end{align*} 
for $1 \leq \alpha,\beta \leq m$. Using the first Bianchi identity, 
we obtain 
\begin{align*} 
&\sum_{i,j,k,l=1}^n (\text{\rm Ric}^{ij} \, g^{kl} - R^{ijkl}) \, \psi_{ik} \, \psi_{jl} \\ 
&= \sum_{\alpha=1}^m \lambda_\alpha^2 \, [\text{\rm Ric}(v_\alpha,v_\alpha) + \text{\rm Ric}(w_\alpha,w_\alpha)] \\ 
&- 2 \sum_{\alpha,\beta=1}^m \lambda_\alpha \, \lambda_\beta \, [R(v_\alpha,v_\beta,w_\alpha,w_\beta) - R(v_\alpha,w_\beta,w_\alpha,v_\beta)] \\ 
&= \sum_{\alpha,\beta=1}^m \lambda_\alpha^2 \, [R(v_\alpha,v_\beta,v_\alpha,v_\beta) + R(v_\alpha,w_\beta,v_\alpha,w_\beta)] \\ 
&+ \sum_{\alpha,\beta=1}^m \lambda_\alpha^2 \, [R(w_\alpha,v_\beta,w_\alpha,v_\beta) + R(w_\alpha,w_\beta,w_\alpha,w_\beta)] \\ 
&- 2 \sum_{\alpha,\beta=1}^m \lambda_\alpha \, \lambda_\beta \, R(v_\alpha,w_\alpha,v_\beta,w_\beta). 
\end{align*} 
This implies 
\begin{align*} 
&\sum_{i,j,k,l=1}^n (\text{\rm Ric}^{ij} \, g^{kl} - R^{ijkl}) \, \psi_{ik} \, \psi_{jl} \\ 
&= \sum_{\alpha \neq \beta} \lambda_\alpha^2 \, [R(v_\alpha,v_\beta,v_\alpha,v_\beta) + R(v_\alpha,w_\beta,v_\alpha,w_\beta)] \\ 
&+ \sum_{\alpha \neq \beta} \lambda_\alpha^2 \, [R(w_\alpha,v_\beta,w_\alpha,v_\beta) + R(w_\alpha,w_\beta,w_\alpha,w_\beta)] \\ 
&- 2 \sum_{\alpha \neq \beta} \lambda_\alpha \, \lambda_\beta \, R(v_\alpha,w_\alpha,v_\beta,w_\beta). 
\end{align*} 
Since $M$ has positive isotropic curvature, we have 
\begin{align*} 
&R(v_\alpha,v_\beta,v_\alpha,v_\beta) + R(v_\alpha,w_\beta,v_\alpha,w_\beta) \\ 
&+ R(w_\alpha,v_\beta,w_\alpha,v_\beta) + R(w_\alpha,w_\beta,w_\alpha,w_\beta) \\ 
&> 2 \, |R(v_\alpha,w_\alpha,v_\beta,w_\beta)|. 
\end{align*} 
for $\alpha \neq \beta$. Since $\sum_{\alpha=1}^m \lambda_\alpha^2 > 0$, it follows that 
\begin{align*} 
&\sum_{i,j,k,l=1}^n (\text{\rm Ric}^{ij} \, g^{kl} - R^{ijkl}) \, \psi_{ik} \, \psi_{jl} \\ 
&> 2 \sum_{\alpha \neq \beta} \lambda_\alpha^2 \, |R(v_\alpha,w_\alpha,v_\beta,w_\beta)| - 2 \sum_{\alpha \neq \beta} |\lambda_\alpha| \, |\lambda_\beta| \, |R(v_\alpha,w_\alpha,v_\beta,w_\beta)| \\
&= \sum_{\alpha \neq \beta} (|\lambda_\alpha| - |\lambda_\beta|)^2 \, |R(v_\alpha,w_\alpha,v_\beta,w_\beta)| \\
&\geq 0
\end{align*} 
at the point $p$. This contradicts (\ref{maximum.point.1}). \\

In odd dimensions, the following result was established by M.~Berger:

\begin{theorem}[M.~Berger \cite{Berger2}]
\label{Betti.number.odd}
Let $M$ be a compact Riemannian manifold of dimension $n \geq 5$. Suppose that $n$ is odd and $M$ has pointwise $\frac{n-3}{4n-9}$-pinched sectional curvatures. Then the second Betti number of $M$ vanishes.
\end{theorem}

\textit{Proof of Theorem \ref{Betti.number.odd}.} 
Suppose that $\psi$ is a non-vanishing harmonic two-form on $M$. The Bochner formula implies that 
\[\Delta \psi_{ik} = \sum_{j=1}^n \text{\rm Ric}_i^j \, \psi_{jk} + \sum_{j=1}^n \text{\rm Ric}_k^j \, \psi_{ij} - 2 \sum_{j,l=1}^n R_{ijkl} \, \psi^{jl}.\] 
As above, we fix a point $p \in M$ where the function $|\psi|^2$ attains its maximum. At the point $p$, we have $|\psi|^2 > 0$ and $\Delta(|\psi|^2) \leq 0$. From this, we deduce that 
\begin{align} 
\label{maximum.point.2}
0 &\geq \Delta \bigg( \sum_{i,k=1}^n \psi_{ik} \, \psi^{ik} \bigg ) \notag \\ 
&\geq 2 \sum_{i,k=1}^n \Delta \psi_{ik} \, \psi^{ik} = 4 \sum_{i,j,k,l=1}^n (\text{\rm Ric}^{ij} \, g^{kl} - R^{ijkl}) \, \psi_{ik} \, \psi_{jl}
\end{align}
at the point $p$. We now write $n = 2m+1$. We can find an orthonormal basis $\{u,v_1,w_1,v_2,w_2,\hdots,v_m,w_m\}$ of $T_p M$ and real numbers $\lambda_1,\hdots,\lambda_m$ such that 
\begin{align*} 
&\psi(u,v_\alpha) = \psi(u,w_\alpha) = 0 \\ 
&\psi(v_\alpha,w_\beta) = \lambda_\alpha \, \delta_{\alpha\beta} \\ 
&\psi(v_\alpha,v_\beta) = \psi(w_\alpha,w_\beta) = 0 
\end{align*}
for $1 \leq \alpha,\beta \leq m$. This implies 
\begin{align*} 
&\sum_{i,j,k,l=1}^n (\text{\rm Ric}^{ij} \, g^{kl} - R^{ijkl}) \, \psi_{ik} \, \psi_{jl} \\ 
&= \sum_{\alpha=1}^m \lambda_\alpha^2 \, [R(u,v_\alpha,u,v_\alpha) + R(u,w_\alpha,u,w_\alpha)] \\ 
&+ \sum_{\alpha \neq \beta} \lambda_\alpha^2 \, [R(v_\alpha,v_\beta,v_\alpha,v_\beta) + R(v_\alpha,w_\beta,v_\alpha,w_\beta)] \\ 
&+ \sum_{\alpha \neq \beta} \lambda_\alpha^2 \, [R(w_\alpha,v_\beta,w_\alpha,v_\beta) + R(w_\alpha,w_\beta,w_\alpha,w_\beta)] \\ 
&- 2 \sum_{\alpha \neq \beta} \lambda_\alpha \, \lambda_\beta \, R(v_\alpha,w_\alpha,v_\beta,w_\beta). 
\end{align*} 
By assumption, $M$ has pointwise $\frac{2m-2}{8m-5}$-pinched sectional curvatures. After rescaling the metric if necessary, we may assume that all sectional curvatures of $M$ at $p$ all lie in the interval $(1,\frac{8m-5}{2m-2}]$. Using Berger's inequality (cf. \cite{Karcher}), we obtain 
\[|R(v_\alpha,w_\alpha,v_\beta,w_\beta)| < \frac{2m-1}{m-1}.\] 
Since $\sum_{\alpha=1}^m \lambda_\alpha^2 > 0$, it follows that 
\begin{align*} 
&\sum_{i,j,k,l=1}^n (\text{\rm Ric}^{ij} \, g^{kl} - R^{ijkl}) \, \psi_{ik} \, \psi_{jl} \\ 
&> (4m-2) \sum_{\alpha=1}^m \lambda_\alpha^2 - \frac{4m-2}{m-1} \sum_{\alpha \neq \beta} |\lambda_\alpha| \, |\lambda_\beta| \\ 
&= \frac{2m-1}{m-1} \sum_{\alpha \neq \beta} (|\lambda_\alpha| - |\lambda_\beta|)^2 \\ 
&\geq 0 
\end{align*} 
at the point $p$. This contradicts (\ref{maximum.point.2}). \\

We note that the pinching constant in Theorem \ref{Betti.number.odd} can be improved for $n = 5$ (see \cite{Berger3}).

\section{The Differentiable Sphere Theorem}

The Topological Sphere Theorem provides a sufficient condition for a Riemannian manifold $M$ to be homeomorphic to $S^n$. We next address the question of whether $M$ is actually diffeomorphic to $S^n$. Various authors have obtained partial results in this direction. The first such result was established in 1966 by D.~Gromoll \cite{Gromoll} and E.~Calabi. Gromoll showed that a simply connected Riemannian manifold whose sectional curvatures lie in the interval $(1,\frac{1}{\delta(n)}]$ is diffeomorphic to $S^n$. The pinching constant $\delta(n)$ depends only on the dimension, and converges to $1$ as $n \to \infty$. In 1971, M.~Sugimoto, K.~Shiohama, and H.~Karcher \cite{Sugimoto-Shiohama} proved the Differentiable Sphere Theorem with a pinching constant $\delta$ independent of $n$ ($\delta = 0.87$). The pinching constant was subsequently improved by E.~Ruh \cite{Ruh1} ($\delta = 0.80$) and by K.~Grove, H.~Karcher, and E.~Ruh \cite{Grove-Karcher-Ruh2} ($\delta = 0.76$). Ruh \cite{Ruh2} proved the Differentiable Sphere Theorem under pointwise pinching assumptions, but with a pinching constant converging to $1$ as $n \to \infty$.

Grove, Karcher, and Ruh \cite{Grove-Karcher-Ruh1},\cite{Grove-Karcher-Ruh2} established an equivariant version of the Differentiable Sphere Theorem, with a pinching constant independent of the dimension ($\delta = 0.98$). The pinching constant was later improved by H.~Im~Hof and E.~Ruh: 

\begin{theorem}[H.~Im~Hof, E.~Ruh \cite{Im-Hof-Ruh}]
There exists a decreasing sequence of real numbers $\delta(n)$ with $\lim_{n \to \infty} \delta(n) = 0.68$ such that the following statement holds: if $M$ is a compact, simply connected $\delta(n)$-pinched Riemannian manifold and $\rho$ is a group homomorphism from a compact Lie group $G$ into the isometry group of $M$, then there exists a diffeomorphism $F: M \to S^n$ and a homomorphism $\sigma: G \to O(n+1)$ such that $F \circ \rho(g) = \sigma(g) \circ F$ for all $g \in G$.
\end{theorem}

In 1982, R.~Hamilton \cite{Hamilton1} introduced fundamental new ideas to this problem. Given a compact Riemannian manifold $(M,g_0)$, Hamilton studied the following evolution equation for the Riemannian metric: 
\begin{equation}
\label{ricci.flow}
\frac{\partial}{\partial t} g(t) = -2 \, \text{\rm Ric}_{g(t)}, \qquad g(0) = g_0. 
\end{equation}
This evolution equation is referred to as the Ricci flow. Hamilton also considered a normalized version of Ricci flow, which differs from the unnormalized flow by a cosmological constant: 
\begin{equation}
\label{normalized.ricci.flow}
\frac{\partial}{\partial t} g(t) = -2 \, \text{\rm Ric}_{g(t)} + \frac{2}{n} \, r_{g(t)} \, g(t), \qquad g(0) = g_0. 
\end{equation}
Here, $r_{g(t)}$ is defined as the mean value of the scalar curvature of $g(t)$. The evolution equations (\ref{ricci.flow}) and (\ref{normalized.ricci.flow}) are essentially equivalent: any solution to equation (\ref{ricci.flow}) can be transformed into a solution of (\ref{normalized.ricci.flow}) by a rescaling procedure (cf. \cite{Hamilton1}).

R.~Hamilton \cite{Hamilton1} proved that the Ricci flow admits a shorttime solution for every initial metric $g_0$ (see also \cite{DeTurck}). Moreover, Hamilton showed that, in dimension $3$, the Ricci flow deforms metrics with positive Ricci curvature to constant curvature metrics:

\begin{theorem}[R.~Hamilton \cite{Hamilton1}]
\label{Hamilton.dim.3}
Let $(M,g_0)$ be a compact three-manifold with positive Ricci curvature. Moreover, let $g(t)$, $t \in [0,T)$, denote the unique maximal solution to the Ricci flow with initial metric $g_0$. Then the rescaled metrics $\frac{1}{4(T-t)} \, g(t)$ converge to a metric of constant sectional curvature $1$ as $t \to T$. In particular, $M$ is diffeomorphic to a spherical space form.
\end{theorem} 

In \cite{Hamilton2}, Hamilton developed powerful techniques for analyzing the global behavior of the Ricci flow. Let $(M,g_0)$ be a compact Riemannian manifold, and let $g(t)$, $t \in [0,T)$, be the unique solution to the Ricci flow with initial metric $g_0$. We denote by $E$ the vector bundle over $M \times (0,T)$ whose fiber over $(p,t) \in M \times (0,T)$ is given by $E_{(p,t)} = T_p M$. The vector bundle admits a natural bundle metric which is defined by $\langle V,W \rangle_h = \langle V,W \rangle_{g(t)}$ for $V,W \in E_{(p,t)}$. Moreover, there is a natural connection $D$ on $E$, which extends the Levi-Civita connection on $TM$. In order to define this connection, we need to specify the covariant time derivative $D_{\frac{\partial}{\partial t}}$. Given two sections $V,W$ of $E$, we define 
\begin{equation} 
\label{covariant.time.derivative}
\langle D_{\frac{\partial}{\partial t}} V,W \rangle_{g(t)} = \langle \frac{\partial}{\partial t} V,W \rangle_{g(t)} - \text{\rm Ric}_{g(t)}(V,W). 
\end{equation}
Note that the connection $D$ is compatible with the bundle metric $h$.

Let $R$ be the curvature tensor of the evolving metric $g(t)$. We may view $R$ as a section of the vector bundle $E^* \otimes E^* \otimes E^* \otimes E^*$. It follows from results of R.~Hamilton \cite{Hamilton2} that $R$ satisfies an evolution equation of the form 
\begin{equation} 
\label{evolution.of.curvature}
D_{\frac{\partial}{\partial t}} R = \Delta R + Q(R). 
\end{equation} 
Here, $D_{\frac{\partial}{\partial t}}$ denotes the covariant time derivative, and $\Delta$ is the Laplacian with respect to the metric $g(t)$. Moreover, $Q(R)$ is defined by 
\begin{equation} 
\label{def.of.Q}
Q(R)_{ijkl} = \sum_{p,q=1}^n R_{ijpq} \, R_{klpq} + 2 \sum_{p,q=1}^n R_{ipkq} \, R_{jplq} - 2 \sum_{p,q=1}^n R_{iplq} \, R_{jpkq}. 
\end{equation} 
Hamilton established a general convergence criterion for the Ricci flow, which reduces the problem to the study of the ODE $\frac{d}{dt} R = Q(R)$ (see \cite{Hamilton2}, Section 5). As an application, Hamilton proved the following convergence theorem in dimension $4$:

\begin{theorem}[R.~Hamilton \cite{Hamilton2}]
\label{Hamilton.dim.4}
Let $(M,g_0)$ be a compact four-manifold with positive curvature operator. Moreover, let $g(t)$, $t \in [0,T)$, denote the unique maximal solution to the Ricci flow with initial metric $g_0$. Then the rescaled metrics $\frac{1}{6(T-t)} \, g(t)$ converge to a metric of constant sectional curvature $1$ as $t \to T$. Consequently, $M$ is diffeomorphic to $S^4$ or $\mathbb{RP}^4$.
\end{theorem} 

H.~Chen \cite{Chen} showed that the conclusion of Theorem \ref{Hamilton.dim.4} holds under the weaker assumption that $(M,g_0)$ has two-positive curvature operator. (That is, the sum of the smallest two eigenvalues of the curvature operator is positive at each point on $M$.) Moreover, Chen proved that any four-manifold with pointwise $1/4$-pinched sectional curvatures has two-positive curvature operator. This implies the following result (see also \cite{Andrews-Nguyen}): 

\begin{theorem}[H.~Chen \cite{Chen}] 
Let $(M,g_0)$ be a compact four-manifold with pointwise $1/4$-pinched sectional curvatures. Let $g(t)$, $t \in [0,T)$, denote the unique maximal solution to the Ricci flow with initial metric $g_0$. Then the rescaled metrics $\frac{1}{6(T-t)} \, g(t)$ converge to a metric of constant sectional curvature $1$ as $t \to T$. 
\end{theorem}

The Ricci flow on manifolds of dimension $n \geq 4$ was first studied by G.~Huisken \cite{Huisken} in 1985 (see also \cite{Margerin1},\cite{Nishikawa}). To describe this result, we decompose the curvature tensor in the usual way as $R_{ijkl} = U_{ijkl} + V_{ijkl} + W_{ijkl}$, where $U_{ijkl}$ denotes the part of the curvature tensor associated with the scalar curvature, $V_{ijkl}$ is the part of the curvature tensor associated with the tracefree Ricci curvature, and $W_{ijkl}$ denotes the Weyl tensor.

\begin{theorem}[G.~Huisken \cite{Huisken}] 
\label{Huiskens.thm}
Let $(M,g_0)$ be a compact Riemannian manifold of dimension $n \geq 4$ with positive scalar curvature. Suppose that the curvature tensor of $(M,g_0)$ satisfies the pointwise pinching condition 
\[|V|^2 + |W|^2 < \delta(n) \, |U|^2,\] 
where $\delta(4) = \frac{1}{5}$, $\delta(5) = \frac{1}{10}$, and 
\[\delta(n) = \frac{2}{(n-2)(n+1)}\] 
for $n \geq 6$. Let $g(t)$, $t \in [0,T)$, denote the unique maximal solution to the Ricci flow with initial metric $g_0$. Then the rescaled metrics $\frac{1}{2(n-1)(T-t)} \, g(t)$ converge to a metric of constant sectional curvature $1$ as $t \to T$. 
\end{theorem}

Note that the curvature condition in Theorem \ref{Huiskens.thm} is preserved by the Ricci flow. Moreover, any manifold $(M,g_0)$ which satisfies the assumptions of Theorem \ref{Huiskens.thm} necessarily has positive curvature operator (see \cite{Huisken}, Corollary 2.5). 

C.~B\"ohm and B.~Wilking proved a convergence result for manifolds with two-positive curvature operator, generalizing Chen's work in dimension $4$:

\begin{theorem}[C.~B\"ohm, B. Wilking \cite{Bohm-Wilking}] 
\label{bw}
Let $(M,g_0)$ is a compact Riemannian manifold with two-positive curvature operator. Let $g(t)$, $t \in [0,T)$, denote the unique maximal solution to the Ricci flow with initial metric $g_0$. Then the rescaled metrics $\frac{1}{2(n-1)(T-t)} \, g(t)$ converge to a metric of constant sectional curvature $1$ as $t \to T$. 
\end{theorem}

C.~Margerin \cite{Margerin2} used the Ricci flow to show that any compact four-manifold which has positive scalar curvature and satisfies the pointwise pinching condition $|W|^2 + |V|^2 < |U|^2$ is diffeomorphic to $S^4$ or $\mathbb{RP}^4$. By combining Margerin's theorem with a conformal deformation of the metric, A.~Chang, M.~Gursky, and P.~Yang were able to replace the pointwise pinching condition by an integral pinching condition. As a result, they obtained a conformally invariant sphere theorem in dimension $4$:

\begin{theorem}[A.~Chang, M.~Gursky, P.~Yang \cite{Chang-Gursky-Yang}] 
\label{chang.gursky.yang.thm}
Let $(M,g_0)$ be a compact four-manifold with positive Yamabe constant. Suppose that $(M,g_0)$ satisfies the integral pinching condition 
\begin{equation} 
\label{integral.pinching.condition}
\int_M (|W|^2 + |V|^2) < \int_M |U|^2.
\end{equation}
Then $M$ is diffeomorphic to $S^4$ or $\mathbb{RP}^4$.
\end{theorem}

Given any compact four-manifold $M$, the Gauss-Bonnet theorem asserts that 
\[\int_M (|U|^2 - |V|^2 + |W|^2) = 32\pi^2 \, \chi(M)\] 
(cf. \cite{Chang-Gursky-Yang}, equation (0.4)). Hence, the condition (\ref{integral.pinching.condition}) is equivalent to 
\begin{equation} 
\int_M |W|^2 < 16\pi^2 \, \chi(M).
\end{equation}
Here, the norm of $W$ is defined by $|W|^2 = \sum_{i,j,k,l=1}^n W_{ijkl} \, W^{ijkl}$.

\section{New invariant curvature conditions for the Ricci flow}

In an important paper \cite{Hamilton8}, R.~Hamilton proved that the Ricci flow preserves positive isotropic curvature in dimension $4$. Moreover, Hamilton studied solutions to the Ricci flow in dimension $4$ with positive isotropic curvature, and analyzed their singularities. Finally, Hamilton \cite{Hamilton8} devised a sophisticated procedure for extending the flow beyond singularities (see also \cite{Chen-Zhu2}, \cite{Perelman1}, \cite{Perelman2}). 

In a recent paper \cite{Brendle-Schoen1}, we proved that positive isotropic curvature is preserved by the Ricci flow in all dimensions. This was shown independently by H.~Nguyen in his doctoral dissertation. Our proof relies on the following algebraic result which is of interest in itself (see \cite{Brendle-Schoen1}, Corollary 10):

\begin{proposition} 
\label{pic.is.preserved.1}
Let $R$ be an algebraic curvature tensor on $\mathbb{R}^n$ with nonnegative isotropic curvature. Moreover, suppose that $\{e_1,e_2,e_3,e_4\}$ is an orthonormal four-frame satisfying 
\[R_{1313} + R_{1414} + R_{2323} + R_{2424} - 2 \, R_{1234} = 0.\] 
Then 
\[Q(R)_{1313} + Q(R)_{1414} + Q(R)_{2323} + Q(R)_{2424} - 2 \, Q(R)_{1234} \geq 0,\] 
where $Q(R)$ is given by (\ref{def.of.Q}).
\end{proposition}

\textit{Sketch of the proof of Proposition \ref{pic.is.preserved.1}.}
Following Hamilton \cite{Hamilton2}, we write $Q(R) = R^2 + R^\#$, where $R^2$ and $R^\#$ are defined by 
\[(R^2)_{ijkl} = \sum_{p,q=1}^n R_{ijpq} \, R_{klpq}\] 
and 
\[(R^\#)_{ijkl} = 2 \sum_{p,q=1}^n R_{ipkq} \, R_{jplq} - 2 \sum_{p,q=1}^n R_{iplq} \, R_{jpkq}.\] 
Note that $R^2$ and $R^\#$ do not satisfy the first Bianchi identity, but $R^2 + R^\#$ does. Since $R$ satisfies the first Bianchi identity, we have 
\begin{align*} 
&(R^\#)_{1313} + (R^\#)_{1414} + (R^\#)_{2323} + (R^\#)_{2424} \\ 
&+ 2 \, (R^\#)_{1342} + 2 \, (R^\#)_{1423} \\ 
&= 2 \sum_{p,q=1}^n (R_{1p1q} + R_{2p2q}) \, (R_{3p3q} + R_{4p4q}) - 2 \sum_{p,q=1}^n R_{12pq} \, R_{34pq} \\ 
&- 2 \sum_{p,q=1}^n (R_{1p3q} + R_{2p4q}) \, (R_{3p1q} + R_{4p2q}) \\ 
&- 2 \sum_{p,q=1}^n (R_{1p4q} - R_{2p3q}) \, (R_{4p1q} - R_{3p2q}). 
\end{align*} 
We claim that the right hand side is nonnegative. To prove this, we define 
\begin{align*} 
I^{(1)} &= \sum_{p,q=1}^4 (R_{1p1q} + R_{2p2q}) \, (R_{3p3q} + R_{4p4q}) - \sum_{p,q=1}^4 R_{12pq} \, R_{34pq} \\ 
&- \sum_{p,q=1}^4 (R_{1p3q} + R_{2p4q}) \, (R_{3p1q} + R_{4p2q}) \\ 
&- \sum_{p,q=1}^4 (R_{1p4q} - R_{2p3q}) \, (R_{4p1q} - R_{3p2q}), 
\end{align*} 
\begin{align*} 
I^{(2)} &= \sum_{p=1}^4 \sum_{q=5}^n (R_{1p1q} + R_{2p2q}) \, (R_{3p3q} + R_{4p4q}) - \sum_{p=1}^4 \sum_{q=5}^n R_{12pq} \, R_{34pq} \\ 
&- \sum_{p=1}^4 \sum_{q=5}^n (R_{1p3q} + R_{2p4q}) \, (R_{3p1q} + R_{4p2q}) \\ 
&- \sum_{p=1}^4 \sum_{q=5}^n (R_{1p4q} - R_{2p3q}) \, (R_{4p1q} - R_{3p2q}), 
\end{align*}
\begin{align*} 
I^{(3)} &= \sum_{p,q=5}^n (R_{1p1q} + R_{2p2q}) \, (R_{3p3q} + R_{4p4q}) - \sum_{p,q=5}^n R_{12pq} \, R_{34pq} \\ 
&- \sum_{p,q=5}^n (R_{1p3q} + R_{2p4q}) \, (R_{3p1q} + R_{4p2q}) \\ 
&- \sum_{p,q=5}^n (R_{1p4q} - R_{2p3q}) \, (R_{4p1q} - R_{3p2q}). 
\end{align*}
We may view the isotropic curvature as a real-valued function on the space of orthonormal four-frames. This function attains its minimum at $\{e_1,e_2,e_3,e_4\}$. Consequently, the first variation at $\{e_1,e_2,e_3,e_4\}$ is zero and the second variation is nonnegative. Using the fact that the first variation is zero, we can show that $I^{(1)} = I^{(2)} = 0$ (see \cite{Brendle-Schoen1}, Propositions 5 and 7). In order to estimate $I^{(3)}$, we consider the following $(n-4) \times (n-4)$ matrices: 
\[\begin{array}{l@{\qquad}l} 
a_{pq} = R_{1p1q} + R_{2p2q}, & b_{pq} = R_{3p3q} +  R_{4p4q}, \\ 
c_{pq} = R_{3p1q} + R_{4p2q}, & d_{pq} = R_{4p1q} - R_{3p2q}, \\ 
e_{pq} = R_{12pq}, & f_{pq} = R_{34pq} 
\end{array}\] 
($5 \leq p,q \leq n$). Since the second variation is nonnegative, the matrix 
\[\begin{bmatrix} B & -F & -C & -D \\ F & B & D & -C \\ -C^T & D^T & A & -E \\ -D^T & -C^T & E & A \end{bmatrix}\] 
is positive semi-definite. From this, we deduce that 
\[I^{(3)} = \text{\rm tr}(AB) + \text{\rm tr}(EF) - \text{\rm tr}(C^2) - \text{\rm tr}(D^2) \geq 0\] 
(see \cite{Brendle-Schoen1}, Proposition 9). Putting these facts together, we conclude that
\begin{align} 
\label{R.sharp.terms}
&(R^\#)_{1313} + (R^\#)_{1414} + (R^\#)_{2323} + (R^\#)_{2424} \notag \\ 
&+ 2 \, (R^\#)_{1342} + 2 \, (R^\#)_{1423} \\ 
&= 2 \, I^{(1)} + 4 \, I^{(2)} + 2 \, I^{(3)} \geq 0. \notag 
\end{align}  
Moreover, we have 
\begin{align} 
\label{R.square.terms}
&(R^2)_{1313} + (R^2)_{1414} + (R^2)_{2323} + (R^2)_{2424} \notag \\ 
&+ 2 \, (R^2)_{1342} + 2 \, (R^2)_{1423} \\ 
&= \sum_{p,q=1}^n (R_{13pq} - R_{24pq})^2 + \sum_{p,q=1}^n (R_{14pq} + R_{23pq})^2 \geq 0 \notag 
\end{align} 
by definition of $R^2$. Adding (\ref{R.sharp.terms}) and (\ref{R.square.terms}), we obtain 
\begin{align} 
&Q(R)_{1313} + Q(R)_{1414} + Q(R)_{2323} + Q(R)_{2424} \notag \\ 
&+ 2 \, Q(R)_{1342} + 2 \, Q(R)_{1423} \geq 0. 
\end{align}  
Since $Q(R)$ satisfies the first Bianchi identity, we conclude that 
\[Q(R)_{1313} + Q(R)_{1414} + Q(R)_{2323} + Q(R)_{2424} - 2 \, Q(R)_{1234} \geq 0,\] 
as claimed. \\

\begin{theorem}[S.~Brendle, R.~Schoen \cite{Brendle-Schoen1}; H.~Nguyen \cite{Nguyen}] 
\label{pic.is.preserved.2}
Let $M$ be a compact manifold of dimension $n \geq 4$, and let $g(t)$, $t \in [0,T)$, be a family of metrics on $M$ evolving under Ricci flow. If $(M,g(0))$ has nonnegative isotropic curvature, then $(M,g(t))$ has nonnegative isotropic curvature for all $t \in [0,T)$.
\end{theorem}

\textit{Proof of Theorem \ref{pic.is.preserved.2}.} 
It follows from Proposition \ref{pic.is.preserved.1} that nonnegative isotropic curvature is preserved by the ODE $\frac{d}{dt} R = Q(R)$. Hence, the assertion follows from Hamilton's maximum principle for systems (see \cite{Hamilton2}). \\

As an application of Proposition \ref{pic.is.preserved.1}, we are able to generalize a theorem of S.~Tachibana \cite{Tachibana}: 

\begin{theorem}[S.~Brendle \cite{Brendle3}]
\label{tachibana.1}
Let $(M,g)$ be a compact Einstein manifold of dimension $n \geq 4$ with positive isotropic curvature. Then $(M,g)$ has constant sectional curvature.
\end{theorem}

\textit{Proof of Theorem \ref{tachibana.1}.} After rescaling the metric if necessary, we may assume that the scalar curvature of $(M,g)$ equals $n(n-1)$. Since $g$ is an Einstein metric, we have $\text{\rm Ric}_{ij} = (n-1) \, g_{ij}$. This implies 
\begin{equation} 
\label{curvature.pde}
\Delta R + Q(R) = 2(n-1) \, R. 
\end{equation}
We define a tensor $S_{ijkl}$ by 
\[S_{ijkl} = R_{ijkl} - \kappa \, (g_{ik} \, g_{jl} - g_{il} \, g_{jk}),\] 
where $\kappa$ is a positive constant. Note that $S$ satisfies all the algebraic properties of the curvature tensor. Let $\kappa$ be the largest constant with the property that $S_{ijkl}$ has nonnegative isotropic curvature. Then there exists a point $p \in M$ and a four-frame $\{e_1,e_2,e_3,e_4\} \subset T_p M$ such that 
\begin{align*} 
&S(e_1,e_3,e_1,e_3) + S(e_1,e_4,e_1,e_4) \\ 
&+ S(e_2,e_3,e_2,e_3) + S(e_2,e_4,e_2,e_4) \\ 
&- 2 \, S(e_1,e_2,e_3,e_4) = 0. 
\end{align*}
Therefore, it follows from Proposition \ref{pic.is.preserved.1} that 
\begin{align} 
\label{Q.term.1}
&Q(S)(e_1,e_3,e_1,e_3) + Q(S)(e_1,e_4,e_1,e_4) \notag \\ 
&+ Q(S)(e_2,e_3,e_2,e_3) + Q(S)(e_2,e_4,e_2,e_4) \\ 
&- 2 \, Q(S)(e_1,e_2,e_3,e_4) \geq 0. \notag 
\end{align}
We next observe that 
\begin{align*} 
Q(S)_{ijkl} 
&= Q(R)_{ijkl} + 2(n-1) \, \kappa^2 \, (g_{ik} \, g_{jl} - g_{il} \, g_{jk}) \\ 
&- 2\kappa \, (\text{\rm Ric}_{ik} \, g_{jl} - \text{\rm Ric}_{il} \, g_{jk} - \text{\rm Ric}_{jk} \, g_{il} + \text{\rm Ric}_{jl} \, g_{ik}), 
\end{align*} 
hence 
\[Q(S)_{ijkl} = Q(R)_{ijkl} + 2(n-1) \, \kappa \, (\kappa - 2) \, (g_{ik} \, g_{jl} - g_{il} \, g_{jk}).\] 
Substituting this into (\ref{Q.term.1}), we obtain 
\begin{align} 
\label{Q.term.2}
&Q(R)(e_1,e_3,e_1,e_3) + Q(R)(e_1,e_4,e_1,e_4) \notag \\ 
&+ Q(R)(e_2,e_3,e_2,e_3) + Q(R)(e_2,e_4,e_2,e_4) \\ 
&- 2 \, Q(R)(e_1,e_2,e_3,e_4) + 8(n-1) \, \kappa \, (\kappa - 2) \geq 0. \notag
\end{align} 
Fix a vector $w \in T_p M$, and consider the geodesic $\gamma(s) = \exp_p(sw)$. Moreover, let $v_j(s)$ be a parallel vector field along $\gamma$ with $v_j(0) = e_j$. The function 
\begin{align*} 
s \mapsto &R(v_1(s),v_3(s),v_1(s),v_3(s)) + R(v_1(s),v_4(s),v_1(s),v_4(s)) \\ 
&+ R(v_2(s),v_3(s),v_2(s),v_3(s)) + R(v_2(s),v_4(s),v_2(s),v_4(s)) \\ 
&- 2 \, R(v_1(s),v_2(s),v_3(s),v_4(s)) - 4\kappa 
\end{align*}
is nonnegative and vanishes for $s = 0$. Hence, the second derivative of that function at $s = 0$ is nonnegative. This implies 
\begin{align*} 
&(D_{w,w}^2 R)(e_1,e_3,e_1,e_3) + (D_{w,w}^2 R)(e_1,e_4,e_1,e_4) \\ 
&+ (D_{w,w}^2 R)(e_2,e_3,e_2,e_3) + (D_{w,w}^2 R)(e_2,e_4,e_2,e_4) \\ 
&- 2 \, (D_{w,w}^2 R)(e_1,e_2,e_3,e_4) \geq 0. 
\end{align*} 
Since $w \in T_p M$ is arbitrary, we conclude that 
\begin{align} 
\label{Laplacian.term}
&(\Delta R)(e_1,e_3,e_1,e_3) + (\Delta R)(e_1,e_4,e_1,e_4) \notag \\ 
&+ (\Delta R)(e_2,e_3,e_2,e_3) + (\Delta R)(e_2,e_4,e_2,e_4) \\ 
&- 2 \, (\Delta R)(e_1,e_2,e_3,e_4) \geq 0. \notag 
\end{align} 
Adding (\ref{Q.term.2}) and (\ref{Laplacian.term}) yields 
\begin{align*}
&R(e_1,e_3,e_1,e_3) + R(e_1,e_4,e_1,e_4) \\ 
&+ R(e_2,e_3,e_2,e_3) + R(e_2,e_4,e_2,e_4) \\ 
&- 2 \, R(e_1,e_2,e_3,e_4) + 4\kappa \, (\kappa - 2) \geq 0. 
\end{align*} 
On the other hand, we have 
\begin{align*}
&R(e_1,e_3,e_1,e_3) + R(e_1,e_4,e_1,e_4) \\ 
&+ R(e_2,e_3,e_2,e_3) + R(e_2,e_4,e_2,e_4) \\ 
&- 2 \, R(e_1,e_2,e_3,e_4) - 4\kappa = 0. 
\end{align*} 
Since $\kappa$ is positive, it follows that $\kappa \geq 1$. Therefore, $S$ has nonnegative isotropic curvature and nonpositive scalar curvature. Hence, Proposition 2.5 in \cite{Micallef-Wang} implies that the Weyl tensor of $(M,g)$ vanishes. \\

In the next step, we apply Theorem \ref{pic.is.preserved.2} to the product manifolds $(M,g(t)) \times \mathbb{R}$ and $(M,g(t)) \times \mathbb{R}^2$.

\begin{theorem}[S.~Brendle, R.~Schoen \cite{Brendle-Schoen1}]
\label{product.with.R}
Let $M$ be a compact manifold of dimension $n \geq 4$, and let $g(t)$, $t \in [0,T)$, be a solution to the Ricci flow on $M$. If $(M,g(0)) \times \mathbb{R}$ has nonnegative isotropic curvature, then $(M,g(t)) \times \mathbb{R}$ has nonnegative isotropic curvature for all $t \in [0,T)$. 
\end{theorem}

\begin{theorem}[S.~Brendle, R.~Schoen \cite{Brendle-Schoen1}]
\label{product.with.R2}
Let $M$ be a compact manifold of dimension $n \geq 4$, and let $g(t)$, $t \in [0,T)$, be a family of metrics on $M$ evolving under Ricci flow. If $(M,g(0)) \times \mathbb{R}^2$ has nonnegative isotropic curvature, then $(M,g(t)) \times \mathbb{R}^2$ has nonnegative isotropic curvature for all $t \in [0,T)$.
\end{theorem}

A similar result holds for products of the form $(M,g(t)) \times S^2(1)$, where $S^2(1)$ denotes a two-dimensional sphere of radius $1$ (see \cite{Brendle1}, Proposition 10).

\begin{theorem}[S.~Brendle \cite{Brendle1}] 
\label{product.with.S2}
Let $M$ be a compact manifold of dimension $n \geq 4$, and let $g(t)$, $t \in [0,T)$, be a solution to the Ricci flow on $M$. If $(M,g(0)) \times S^2(1)$ has nonnegative isotropic curvature, then $(M,g(t)) \times S^2(1)$ has nonnegative isotropic curvature for all $t \in [0,T)$. 
\end{theorem}

Theorem \ref{product.with.S2} is quite subtle, as the manifolds $(M,g(t)) \times S^2(1)$ do not form a solution to the Ricci flow. 

Theorems \ref{product.with.R} -- \ref{product.with.S2} provide us with various curvature conditions that are preserved by the Ricci flow. We now discuss these curvature conditions in more detail. Let $M$ be a Riemannian manifold of dimension $n \geq 4$. The product $M \times \mathbb{R}$ has nonnegative isotropic curvature if and only if 
\[R_{1313} + \lambda^2 \, R_{1414} + R_{2323} + \lambda^2 \, R_{2424} - 2\lambda \, R_{1234} \geq 0\] 
for all orthonormal four-frames $\{e_1,e_2,e_3,e_4\}$ and all $\lambda \in [-1,1]$ (see \cite{Brendle1}, Proposition 4). 
Similarly, the product $M \times \mathbb{R}^2$ has nonnegative isotropic curvature if and only if 
\[R_{1313} + \lambda^2 \, R_{1414} + \mu^2 \, R_{2323} + \lambda^2\mu^2 \, R_{2424} - 2\lambda\mu \, R_{1234} \geq 0\] 
for all orthonormal four-frames $\{e_1,e_2,e_3,e_4\}$ and all $\lambda,\mu \in [-1,1]$ (see \cite{Brendle-Schoen1}, Proposition 21). Finally, the product $M \times S^2(1)$ has nonnegative isotropic curvature if and only if 
\begin{align*} 
&R_{1313} + \lambda^2 \, R_{1414} + \mu^2 \, R_{2323} + \lambda^2\mu^2 \, R_{2424} \\ 
&- 2\lambda\mu \, R_{1234} + (1 - \lambda^2) \, (1 - \mu^2) \geq 0 
\end{align*} 
for all orthonormal four-frames $\{e_1,e_2,e_3,e_4\}$ and all $\lambda,\mu \in [-1,1]$ (cf. \cite{Brendle1}, Proposition 7).

We can also characterize these curvature conditions using complex notation. The product $M \times \mathbb{R}$ has nonnegative isotropic curvature if and only if $R(z,w,\bar{z},\bar{w}) \geq 0$ for all vectors $z,w \in T_p^{\mathbb{C}} M$ satisfying $g(z,z) \, g(w,w) - g(z,w)^2 = 0$. Moreover, the product $M \times \mathbb{R}^2$ has nonnegative isotropic curvature if and only if $R(z,w,\bar{z},\bar{w}) \geq 0$ for all vectors $z,w \in T_p^{\mathbb{C}} M$ (see \cite{Micallef-Wang}, Remark 3.3).

Combining these results with earlier work of Hamilton \cite{Hamilton2} and of B\"ohm and Wilking \cite{Bohm-Wilking}, we obtain the following theorem:

\begin{theorem}[S.~Brendle, R.~Schoen \cite{Brendle-Schoen1}]
\label{convergence.1}
Let $(M,g_0)$ be a compact Riemannian manifold of dimension $n \geq 4$ such that 
\begin{equation} 
\label{curvature.condition.1}
R_{1313} + \lambda^2 \, R_{1414} + \mu^2 \, R_{2323} + \lambda^2\mu^2 \, R_{2424} - 2\lambda\mu \, R_{1234} > 0 
\end{equation} 
for all orthonormal four-frames $\{e_1,e_2,e_3,e_4\}$ and all $\lambda,\mu \in [-1,1]$. Let $g(t)$, $t \in [0,T)$, denote the unique maximal solution to the Ricci flow with initial metric $g_0$. Then the rescaled metrics $\frac{1}{2(n-1)(T-t)} \, g(t)$ converge to a metric of constant sectional curvature $1$ as $t \to T$. 
\end{theorem}

It follows from Berger's inequality that every manifold with pointwise $1/4$-pinched sectional curvatures satisfies (\ref{curvature.condition.1}). Hence, we can draw the following conclusion:

\begin{corollary}[S.~Brendle, R.~Schoen \cite{Brendle-Schoen1}] 
\label{diffeo.sphere.thm}
Let $(M,g_0)$ be a compact Riemannian manifold with pointwise $1/4$-pinched sectional curvatures. Moreover, let $g(t)$, $t \in [0,T)$, denote the unique maximal solution to the Ricci flow with initial metric $g_0$. Then the rescaled metrics $\frac{1}{2(n-1)(T-t)} \, g(t)$ converge to a metric of constant sectional curvature $1$ as $t \to T$. 
\end{corollary}

Both Theorem \ref{bw} and Theorem \ref{convergence.1} are subcases of a more general convergence theorem for the Ricci flow: 

\begin{theorem}[S.~Brendle \cite{Brendle1}]
\label{convergence.2}
Let $(M,g_0)$ be a compact Riemannian manifold of dimension $n \geq 4$ such that 
\begin{equation} 
\label{curvature.condition.2}
R_{1313} + \lambda^2 \, R_{1414} + R_{2323} + \lambda^2 \, R_{2424} - 2\lambda \, R_{1234} > 0 
\end{equation} 
for all orthonormal four-frames $\{e_1,e_2,e_3,e_4\}$ and all $\lambda \in [-1,1]$. Let $g(t)$, $t \in [0,T)$, denote the unique maximal solution to the Ricci flow with initial metric $g_0$. Then the rescaled metrics $\frac{1}{2(n-1)(T-t)} \, g(t)$ converge to a metric of constant sectional curvature $1$ as $t \to T$.
\end{theorem}

To conclude this section, we provide a diagram showing the logical implications 
among the following curvature conditions: 
\begin{itemize}
\item[(C1)] $M$ has $1/4$-pinched sectional curvatures 
\item[(C2)] $M$ has nonnegative sectional curvature 
\item[(C3)] $M$ has two-nonnegative flag curvature; that is, $R_{1313} + R_{2323} \geq 0$ for all orthonormal three-frames $\{e_1,e_2,e_3\}$
\item[(C4)] $M$ has nonnegative scalar curvature
\item[(C5)] $M \times \mathbb{R}^2$ has nonnegative isotropic curvature 
\item[(C6)] $M \times S^2(1)$ has nonnegative isotropic curvature 
\item[(C7)] $M \times \mathbb{R}$ has nonnegative isotropic curvature 
\item[(C8)] $M$ has nonnegative isotropic curvature 
\item[(C9)] $M$ has nonnegative curvature operator 
\item[(C10)] $M$ has two-nonnegative curvature operator 
\end{itemize}
Note that conditions (C4) -- (C10) are preserved by the Ricci flow, but (C1) -- (C3) are not.

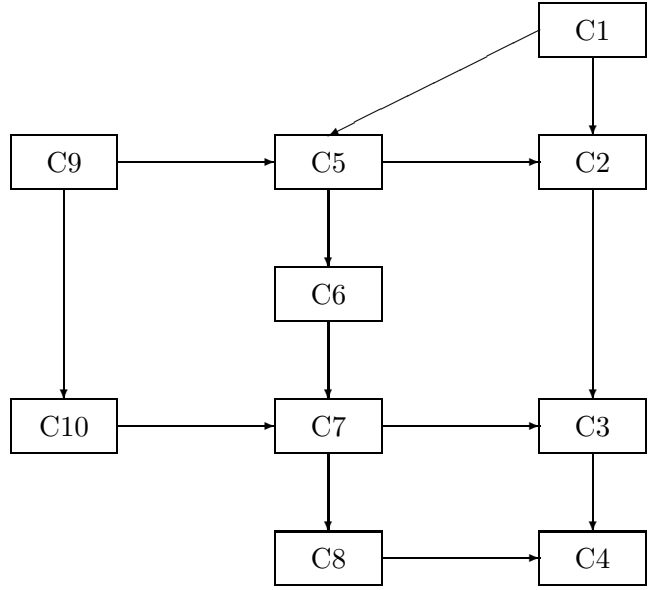
\begin{center}
\begin{picture}(270,250)
\put(0,160){\framebox(40,20)[c]{C9}}
\put(0,60){\framebox(40,20)[c]{C10}}
\put(100,160){\framebox(40,20)[c]{C5}}
\put(100,110){\framebox(40,20)[c]{C6}}
\put(100,60){\framebox(40,20)[c]{C7}}
\put(100,10){\framebox(40,20)[c]{C8}}
\put(200,210){\framebox(40,20)[c]{C1}}
\put(200,160){\framebox(40,20)[c]{C2}}
\put(200,60){\framebox(40,20)[c]{C3}}
\put(200,10){\framebox(40,20)[c]{C4}}
\put(20,160){\vector(0,-1){80}}
\put(120,160){\vector(0,-1){30}}
\put(120,110){\vector(0,-1){30}}
\put(120,60){\vector(0,-1){30}}
\put(200,220){\vector(-2,-1){80}}
\put(220,210){\vector(0,-1){30}}
\put(220,160){\vector(0,-1){80}}
\put(220,60){\vector(0,-1){30}}
\put(40,170){\vector(1,0){60}}
\put(140,170){\vector(1,0){60}}
\put(40,70){\vector(1,0){60}}
\put(140,70){\vector(1,0){60}}
\put(140,20){\vector(1,0){60}}
\end{picture}
\end{center}

\section{Rigidity results and the classification of weakly $1/4$-pinched manifolds}

In this section, we describe various rigidity results. The following theorem is based on the strict maximum principle and plays a key role in our analysis:

\begin{theorem}[S.~Brendle, R.~Schoen \cite{Brendle-Schoen2}] 
\label{zero.isotropic.curvature.1}
Let $M$ be a compact manifold of dimension $n \geq 4$, and let $g(t)$, $t \in [0,T)$, be a solution to the Ricci flow on $M$ with nonnegative isotropic curvature. Moreover, we fix a time $\tau \in (0,T)$. Then the set of all four-frames $\{e_1,e_2,e_3,e_4\}$ that are orthonormal with respect to $g(\tau)$ and satisfy 
\begin{align*} 
&R_{g(\tau)}(e_1,e_3,e_1,e_3) + R_{g(\tau)}(e_1,e_4,e_1,e_4) \\ 
&+ R_{g(\tau)}(e_2,e_3,e_2,e_3) + R_{g(\tau)}(e_2,e_4,e_2,e_4) \\ 
&- 2 \, R_{g(\tau)}(e_1,e_2,e_3,e_4) = 0 
\end{align*} 
is invariant under parallel transport.
\end{theorem}

\textit{Sketch of the proof of Theorem \ref{zero.isotropic.curvature.1}.} 
Let $E$ denote the vector bundle defined in Section 3, and let $P$ be the orthonormal frame bundle of $E$; that is, the fiber of $P$ over a point $(p,t) \in M \times (0,T)$ consists of all $n$-frames $\{e_1, \hdots, e_n\} \subset T_p M$ that are orthonormal with respect to the metric $g(t)$. Note that $P$ is a principal $O(n)$-bundle over $M \times (0,T)$. Let $\pi$ denote the projection from $P$ to $M \times (0,T)$. For each $t \in (0,T)$, we denote by $P_t = \pi^{-1}(M \times \{t\})$ the time $t$ slice of $P$. 

The connection $D$ defines a horizontal distribution on $P$. For each point $\underline{e} = \{e_1, \hdots, e_n\} \in P$, the tangent space $T_{\underline{e}} P$ splits as a direct sum $T_{\underline{e}} P = \mathbb{H}_{\underline{e}} \oplus \mathbb{V}_{\underline{e}}$, where $\mathbb{H}_{\underline{e}}$ and $\mathbb{V}_{\underline{e}}$ denote the horizontal and 
vertical subspaces at $\underline{e}$, respectively. We next define a collection of horizontal vector fields $\tilde{X}_1, \hdots, \tilde{X}_n, \tilde{Y}$ on $P$. For each $j = 1, \hdots, n$, the value of $\tilde{X}_j$ at a point $\underline{e} = \{e_1, \hdots, e_n\} \in P$ is given by the horizontal lift of the vector $e_j$. Similarly, the value of $\tilde{Y}$ at a point $\underline{e} = \{e_1, \hdots, e_n\} \in P$ is given by the horizontal lift of the vector $\frac{\partial}{\partial t}$. Note that the vector fields $\tilde{X}_1, \hdots, \tilde{X}_n$ are tangential to $P_t$. 

We define a function $u: P \to \mathbb{R}$ by 
\begin{align*} 
u: \underline{e} = \{e_1,\hdots,e_n\} \mapsto \; &R(e_1,e_3,e_1,e_3) + R(e_1,e_4,e_1,e_4) \\ 
&+ R(e_2,e_3,e_2,e_3) + R(e_2,e_4,e_2,e_4) \\ &- 2 \, R(e_1,e_2,e_3,e_4), 
\end{align*} 
where $R$ denotes the Riemann curvature tensor of the evolving metric $g(t)$. By assumption, the function $u: P \to \mathbb{R}$ is nonnegative. Using (\ref{evolution.of.curvature}), we obtain 
\begin{align*} 
\tilde{Y}(u) - \sum_{j=1}^n \tilde{X}_j(\tilde{X}_j(u)) &= Q(R)(e_1,e_3,e_1,e_3) + Q(R)(e_1,e_4,e_1,e_4) \\ 
&+ Q(R)(e_2,e_3,e_2,e_3) + Q(R)(e_2,e_4,e_2,e_4) \\ 
&- 2 \, Q(R)(e_1,e_2,e_3,e_4) 
\end{align*} 
(see \cite{Brendle-Schoen2}, Lemma 6). Moreover, there exists a positive constant $K$ such that 
\begin{align*} 
&Q(R)(e_1,e_3,e_1,e_3) + Q(R)(e_1,e_4,e_1,e_4) \\ 
&+ Q(R)(e_2,e_3,e_2,e_3) + Q(R)(e_2,e_4,e_2,e_4) \\ 
&- 2 \, Q(R)(e_1,e_2,e_3,e_4) \\ 
&\geq K \, \inf_{\xi \in \mathbb{V}_{\underline{e}}, \, |\xi| \leq 1} (D^2 u)(\xi,\xi) - K \, \sup_{\xi \in \mathbb{V}_{\underline{e}}, \, |\xi| \leq 1} Du(\xi) - K \, u 
\end{align*} 
(see \cite{Brendle-Schoen2}, Lemma 7). Putting these facts together, we obtain 
\begin{align*} 
&\tilde{Y}(u) - \sum_{j=1}^n \tilde{X}_j(\tilde{X}_j(u)) \\ 
&\geq K \, \inf_{\xi \in \mathbb{V}_{\underline{e}}, \, |\xi| \leq 1} (D^2 u)(\xi,\xi) - K \, \sup_{\xi \in \mathbb{V}_{\underline{e}}, \, |\xi| \leq 1} Du(\xi) - K \, u. 
\end{align*}
Since $u$ satisfies this inequality, we may apply a variant of Bony's strict maximum principle for degenerate elliptic equations (cf. \cite{Bony}). Hence, if $\tilde{\gamma}: [0,1] \to P_\tau$ is a horizontal curve such that $\tilde{\gamma}(0)$ lies in the zero set of the function $u$, then $\tilde{\gamma}(1)$ lies in the zero set of the function $u$ (see \cite{Brendle-Schoen2}, Proposition 5). From this, the assertion follows. \\

If we apply Proposition \ref{zero.isotropic.curvature.1} to the product manifolds $(M,g(t)) \times S^1$, then we obtain the following result:

\begin{corollary}[S.~Brendle, R.~Schoen \cite{Brendle-Schoen2}] 
\label{zero.isotropic.curvature.2}
Let $M$ be a compact manifold of dimension $n \geq 4$. Moreover, let $g(t)$, $t \in [0,T)$, be a 
solution to the Ricci flow on $M$ with the property that $(M,g(t)) \times \mathbb{R}$ has nonnegative isotropic curvature. 
Fix real numbers $\tau \in (0,T)$ and $\lambda \in [-1,1]$. Then the set of all four-frames $\{e_1,e_2,e_3,e_4\}$ that are orthonormal with respect to $g(\tau)$ and satisfy 
\begin{align*} 
&R_{g(\tau)}(e_1,e_3,e_1,e_3) + \lambda^2 \, R_{g(\tau)}(e_1,e_4,e_1,e_4) \\ 
&+ R_{g(\tau)}(e_2,e_3,e_2,e_3) + \lambda^2 \, R_{g(\tau)}(e_2,e_4,e_2,e_4) \\ 
&- 2\lambda \, R_{g(\tau)}(e_1,e_2,e_3,e_4) = 0 
\end{align*} 
is invariant under parallel transport.
\end{corollary}

Theorem \ref{zero.isotropic.curvature.1} and Corollary \ref{zero.isotropic.curvature.2} can be used to prove various rigidity results. For example, we can extend Theorem \ref{tachibana.1} as follows:

\begin{theorem}[S.~Brendle \cite{Brendle3}]
\label{tachibana.2}
Let $(M,g)$ be a compact Einstein manifold of dimension $n \geq 4$ with nonnegative isotropic curvature. Then $(M,g)$ is locally symmetric.
\end{theorem}

In the next step, we classify all Riemannian manifolds $(M,g_0)$ with the property that $(M,g_0) \times \mathbb{R}$ has nonnegative isotropic curvature: 

\begin{theorem} 
\label{borderline.case.1}
Let $(M,g_0)$ be a compact, locally irreducible Riemannian manifold of dimension $n \geq 4$. Suppose that $(M,g_0) \times \mathbb{R}$ has nonnegative isotropic curvature. Moreover, let $g(t)$, $t \in [0,T)$, denote the unique maximal solution to the Ricci flow with initial metric $g_0$. Then one of the following statements holds: 
\begin{itemize}
\item[(i)] The rescaled metrics $\frac{1}{2(n-1)(T-t)} \, g(t)$ converge to a metric of constant sectional curvature $1$ as $t \to T$.
\item[(ii)] $n = 2m$ and the universal cover of $(M,g_0)$ is a K\"ahler manifold.
\item[(iii)] $(M,g_0)$ is locally symmetric.
\end{itemize}
\end{theorem}

\textit{Sketch of the proof of Theorem \ref{borderline.case.1}.} 
By assumption, the manifold $(M,g_0)$ is locally irreducible and has nonnegative Ricci curvature. By a theorem of Cheeger and Gromoll, the universal cover of $M$ is compact (see \cite{Cheeger-Gromoll} or \cite{Petersen}, p.~288).

If $(M,g_0)$ is locally symmetric, we are done. Hence, we will assume that $(M,g_0)$ is not locally symmetric. By continuity, there exists a real number $\delta \in (0,T)$ such that $(M,g(t))$ is locally irreducible and non-symmetric for all $t \in (0,\delta)$. By Berger's holonomy theorem, there are three possibilities:

\textit{Case 1:} Suppose that $\text{\rm Hol}^0(M,g(\tau)) = SO(n)$ for some $\tau \in (0,\delta)$. In this case, it follows from Corollary \ref{zero.isotropic.curvature.2} that 
\begin{align*} 
&R_{g(\tau)}(e_1,e_3,e_1,e_3) + \lambda^2 \, R_{g(\tau)}(e_1,e_4,e_1,e_4) \\ 
&+ R_{g(\tau)}(e_2,e_3,e_2,e_3) + \lambda^2 \, R_{g(\tau)}(e_2,e_4,e_2,e_4) \\ 
&- 2\lambda \, R_{g(\tau)}(e_1,e_2,e_3,e_4) > 0 
\end{align*} 
for all orthonormal four-frames $\{e_1,e_2,e_3,e_4\}$ and all $\lambda \in [-1,1]$. By Theorem \ref{convergence.2}, the 
rescaled metrics $\frac{1}{2(n-1)(T-t)} \, g(t)$ converge to a metric of constant sectional curvature $1$ as $t \to T$.

\textit{Case 2:} Suppose that $n = 2m$ and $\text{\rm Hol}^0(M,g(t)) = U(m)$ for all $t \in (0,\delta)$. In this case, the universal cover of $(M,g(t))$ is a K\"ahler manifold for all $t \in (0,\delta)$. Since $g(t) \to g_0$ in $C^\infty$, we conclude that the universal cover of $(M,g_0)$ is a K\"ahler manifold.

\textit{Case 3:} Suppose that $n = 4m \geq 8$ and $\text{\rm Hol}^0(M,g(\tau)) = \text{\rm Sp}(m) \cdot \text{\rm Sp}(1)$ for some $\tau \in (0,\delta)$. In this case, the universal cover of $(M,g(\tau))$ is a compact quaternionic-K\"ahler manifold. In particular, $(M,g(\tau))$ is an Einstein manifold. Since $(M,g(\tau))$ has nonnegative isotropic curvature, Theorem \ref{tachibana.2} implies that $(M,g(\tau))$ is locally symmetric. This is a contradiction. \\

In the special case that $(M,g_0)$ has weakly $1/4$-pinched sectional curvatures, we can draw the following conclusion:

\begin{corollary}[S.~Brendle, R.~Schoen \cite{Brendle-Schoen2}] 
\label{borderline.case.2}
Assume that $(M,g_0)$ has weakly $1/4$-pinched sectional curvatures in the sense that $0 \leq K(\pi_1) \leq 4 \, K(\pi_2)$ for all points $p \in M$ and all two-planes $\pi_1,\pi_2 \subset T_p M$. Moreover, we assume that $(M,g_0)$ is not locally symmetric. Finally, let $g(t)$, $t \in [0,T)$, denote the unique maximal solution to the Ricci flow with initial metric $g_0$. Then the rescaled metrics $\frac{1}{2(n-1)(T-t)} \, g(t)$ converge to a metric of constant sectional curvature $1$ as $t \to T$.
\end{corollary}

Finally, we briefly discuss the problem of classifying manifolds with almost $1/4$-pinched sectional curvature. This question was first studied by M.~Berger \cite{Berger4}. Berger showed that for each even integer $n$ there exists a positive real number $\varepsilon(n)$ with the following property: if $M$ is a compact, simply connected Riemannian manifold of dimension $n$ whose sectional curvatures lie in the interval $(1,4+\varepsilon(n)]$, then $M$ is homeomorphic to $S^n$ or diffeomorphic to a compact symmetric space of rank one. 

Building upon earlier work of J.P.~Bourguignon \cite{Bourguignon}, W.~Seaman \cite{Seaman} proved that a compact, simply connected four-manifold whose sectional curvatures lie in the interval $(0.188,1]$ is homeomorphic to $S^4$ or $\mathbb{CP}^2$. U.~Abresch and W.~Meyer \cite{Abresch-Meyer} showed that a compact, simply connected, odd-dimensional Riemannian manifold whose sectional curvatures lie in the interval $(1,4(1+10^{-6})^2]$ is homeomorphic to a sphere. 

Using Theorem \ref{borderline.case.1} and Cheeger-Gromov compactness theory, P.~Petersen and T.~Tao \cite{Petersen-Tao} proved that any compact, simply connected Riemannian manifold of dimension $n$ whose sectional curvatures lie in the interval $(1,4+\varepsilon(n)]$ is diffeomorphic to a sphere or a compact symmetric space of rank one. Here, $\varepsilon(n)$ is a positive real number which depends only on $n$.

\section{Hamilton's differential Harnack inequality for the Ricci flow} 

In 1993, R.~Hamilton \cite{Hamilton4} established a differential Harnack inequality for solutions to the Ricci flow with nonnegative curvature operator (see \cite{Hamilton3} for an earlier result in dimension $2$). In this section, we describe this inequality, as well as some applications. Let $(M,g(t))$, $t \in (0,T)$, be a family of complete Riemannian manifolds evolving under Ricci flow. Following R.~Hamilton \cite{Hamilton1}, we define 
\[P_{ijk} = D_i \text{\rm Ric}_{jk} - D_j \text{\rm Ric}_{ik}\] 
and 
\[M_{ij} = \Delta \text{\rm Ric}_{ij} - \frac{1}{2} \, D_{i,j}^2 \text{\rm scal} + 2 \, R_{ikjl}  \, \text{\rm Ric}^{kl} - \text{\rm Ric}_i^k \, \text{\rm Ric}_{jk} + \frac{1}{2t} \, \text{\rm Ric}_{ij}.\] 
Hamilton's matrix Harnack inequality states:

\begin{theorem}[R.~Hamilton \cite{Hamilton4}]
\label{Hamiltons.Harnack.inequality}
Let $(M,g(t))$, $t \in (0,T)$, be a solution to the Ricci flow with uniformly bounded curvature and nonnegative curvature operator. Then 
\[M(w,w) + 2 \, P(v,w,w) + R(v,w,v,w) \geq 0\] 
for all points $(p,t) \in M \times (0,T)$ and all vectors $v,w \in T_p M$. 
\end{theorem}

Taking the trace over $w$, Hamilton obtained a gradient estimate for the scalar curvature:

\begin{corollary}[R.~Hamilton \cite{Hamilton4}] 
\label{trace.Harnack.inequality}
Assume that $(M,g(t))$, $t \in (0,T)$, is a solution to the Ricci flow with uniformly bounded curvature and nonnegative curvature operator. Then 
\[\frac{\partial}{\partial t} \text{\rm scal} + \frac{1}{t} \, \text{\rm scal} + 2 \, \partial_i \text{\rm scal} \, v^i + 2 \, \text{\rm Ric}(v,v) \geq 0\] 
for all points $(p,t) \in M \times (0,T)$ and all vectors $v \in T_p M$. 
\end{corollary}

We note that H.D.~Cao \cite{Cao} has established a differential Harnack inequality for solutions to the K\"ahler-Ricci flow with nonnegative holomorphic bisectional curvature. In \cite{Brendle2}, it was shown that Hamilton's Harnack inequality holds under the weaker assumption that $(M,g(t)) \times \mathbb{R}^2$ has nonnegative isotropic curvature:

\begin{theorem}[S.~Brendle \cite{Brendle2}]
\label{generalization.of.Hamiltons.Harnack.inequality}
Let $(M,g(t))$, $t \in (0,T)$, be a solution to the Ricci flow with uniformly bounded curvature. Moreover, suppose that the product $(M,g(t)) \times \mathbb{R}^2$ has nonnegative isotropic curvature for all $t \in (0,T)$. Then 
\[M(w,w) + 2 \, P(v,w,w) + R(v,w,v,w) \geq 0\] 
for all points $(p,t) \in M \times (0,T)$ and all vectors $v,w \in T_p M$. 
\end{theorem}

The Harnack inequality has various applications. For example, it can be used to show that any Type II singularity model with nonnegative curvature operator and strictly positive Ricci curvature must be a steady Ricci soliton (see \cite{Hamilton5}). A Riemannian manifold $(M,g)$ is called a gradient Ricci soliton if there exists a smooth function $f: M \to \mathbb{R}$ and a constant $\rho$ such that $\text{\rm Ric}_{ij} = \rho \, g_{ij} + D_{i,j}^2 f$. Depending on the sign of $\rho$, a gradient Ricci soliton is called shrinking ($\rho > 0$), steady ($\rho = 0$), or expanding ($\rho < 0$). As in Theorem \ref{generalization.of.Hamiltons.Harnack.inequality}, we can replace the condition that $M$ has nonnegative curvature operator by the weaker condition that $M \times \mathbb{R}^2$ has nonnegative isotropic curvature:

\begin{proposition}[S.~Brendle \cite{Brendle2}]
\label{singularity.model.of.type.2}
Let $(M,g(t))$, $t \in (-\infty,T)$, be a solution to the Ricci flow which is complete and simply connected. We assume that $(M,g(t)) \times \mathbb{R}^2$ has nonnegative isotropic curvature and $(M,g(t))$ has positive Ricci curvature. Moreover, suppose that there exists a point $(p_0,t_0) \in M \times (-\infty,T)$ such that 
\[\text{\rm scal}_{g(t)}(p) \leq \text{\rm scal}_{g(t_0)}(p_0)\] 
for all points $(p,t) \in M \times (-\infty,T)$. Then $(M,g(t_0))$ is a steady gradient Ricci soliton.
\end{proposition} 

\begin{proposition}[S.~Brendle \cite{Brendle2}]
\label{singularity.model.of.type.3}
Let $(M,g(t))$, $t \in (0,T)$, be a solution to the Ricci flow which is complete and simply connected. We assume that $(M,g(t)) \times \mathbb{R}^2$ has nonnegative isotropic curvature and $(M,g(t))$ has positive Ricci curvature. Moreover, suppose that there exists a point $(p_0,t_0) \in M \times (0,T)$ such that 
\[t \cdot \text{\rm scal}_{g(t)}(p) \leq t_0 \cdot \text{\rm scal}_{g(t_0)}(p_0)\] 
for all points $(p,t) \in M \times (0,T)$. Then $(M,g(t_0))$ is an expanding gradient Ricci soliton.
\end{proposition}

\section{Compactness of pointwise pinched manifolds}

R.~Hamilton \cite{Hamilton6} has shown that a convex hypersurface with pinched second fundamental form is necessarily compact. B.~Chen and X.~Zhu proved an intrinsic analogue of this result. More precisely, they showed that a complete Riemannian manifold which satisfies a suitable pointwise pinching condition is compact: 

\begin{theorem}[B.~Chen, X.~Zhu \cite{Chen-Zhu1}]
\label{chen.zhu.compactness}
Let $(M,g_0)$ be a complete Riemannian manifold of dimension $n \geq 4$. Assume that the scalar curvature of $(M,g_0)$ is uniformly bounded and positive. Moreover, suppose that $(M,g_0)$ satisfies the pointwise pinching condition 
\[|V|^2 + |W|^2 < \delta(n) \, (1 - \varepsilon)^2 \, |U|^2,\] 
where $\varepsilon$ is a positive real number and $\delta(n)$ denotes the pinching constant defined in Theorem \ref{Huiskens.thm}. Then $M$ is compact.
\end{theorem}

Theorem \ref{chen.zhu.compactness} was generalized by L.~Ni and B.~Wu (see \cite{Ni-Wu}, Theorem 3.1). In the remainder of this section, we prove another generalization of Theorem \ref{chen.zhu.compactness}. To that end, we need two results due to L.~Ni \cite{Ni} and L.~Ma and D.~Chen \cite{Ma-Chen}: 

\begin{proposition}[L.~Ni \cite{Ni}]
\label{curvature.decay.on.steady.solitons}
Let $(M,g)$ be a steady gradient Ricci soliton which is complete and non-compact. Suppose that there exists a point $p_0 \in M$ such that $0 < \text{\rm scal}(p) \leq \text{\rm scal}(p_0)$ for all points $p \in M$. Moreover, we assume that $\text{\rm Ric} \geq \varepsilon \, \text{\rm scal} \: g$ for some constant $\varepsilon > 0$. Then there exists a constant $\alpha > 0$ such that 
\[\text{\rm scal}(p) \leq e^{-\alpha \, d(p_0,p)} \, \text{\rm scal}(p_0)\] 
for $d(p_0,p) \geq 1$.
\end{proposition}

\textit{Proof of Proposition \ref{curvature.decay.on.steady.solitons}.} 
Since $(M,g)$ is a steady gradient Ricci soliton, there exists a smooth function $f: M \to \mathbb{R}$ such that $\text{\rm Ric}_{ij} = D_{i,j}^2 f$. This implies 
\begin{align*} 
0 &= \partial_i \text{\rm scal} - 2 \, g^{kl} \, D_i \text{\rm Ric}_{kl} + 2 \, g^{kl} \, D_k \text{\rm Ric}_{il} \\ 
&= \partial_i \text{\rm scal} - 2 \, g^{kl} \, D_{i,k,l}^3 f + 2 \, g^{kl} \, D_{k,i,l}^3 f \\ 
&= \partial_i \text{\rm scal} + 2 \, g^{kl} \, R_{ikjl} \, \partial^j f \\ 
&= \partial_i \text{\rm scal} + 2 \, \text{\rm Ric}_{ij} \, \partial^j f. 
\end{align*} 
By assumption, the scalar curvature attains its maximum at the point $p_0$. This implies $\partial_i \text{\rm scal}(p_0) = 0$. Since $M$ has positive Ricci curvature, it follows that $\partial_i f(p_0) = 0$. Hence, the point $p_0$ is a critical point of the function $f$. Since $f$ is strictly convex, we conclude that $f$ has no critical points other than $p_0$. Let $\sigma$ be a positive real number such that $\text{\rm Ric} \geq \sigma \, g$ for $d(p_0,p) \leq 1$. This implies $f(p) - f(p_0) \geq \frac{1}{2} \, \sigma \, d(p_0,p)^2$ for $d(p_0,p) \leq 1$. Moreover, we have $f(p) - f(p_0) \geq \frac{1}{2} \, \sigma \, d(p_0,p)$ for $d(p_0,p) \geq 1$.

Using the inequality $\text{\rm Ric} \geq \varepsilon \, \text{\rm scal} \: g$, we obtain 
\begin{align*} 
\partial_i(e^{2\varepsilon \, f} \, \text{\rm scal}) \, \partial^i f 
&= e^{2\varepsilon f} \, (\partial_i \text{\rm scal} \, \partial^i f + 2\varepsilon \, \text{\rm scal} \: |Df|^2) \\ 
&\leq e^{2\varepsilon f} \, (\partial_i \text{\rm scal} \, \partial^i f + 2 \, \text{\rm Ric}_{ij} \, \partial^i f \, \partial^j f) \\ 
&= 0. 
\end{align*} 
Fix a point $p \in M$, and let $\gamma: [0,\infty) \to M$ be the solution of the ODE $\gamma'(s) = -\partial^i f(\gamma(s)) \, \frac{\partial}{\partial x^i}$ with initial condition $\gamma(0) = p$. It is easy to see that $\gamma(s) \to p_0$ as $s \to \infty$. Moreover, the function $s \mapsto e^{2\varepsilon \, f(\gamma(s))} \, \text{\rm scal}(\gamma(s))$ is monotone increasing. Thus, we conclude that 
\[\text{\rm scal}(p) \leq e^{-2\varepsilon (f(p)-f(p_0))} \, \text{\rm scal}(p_0)\] 
for all points $p \in M$. This implies 
\[\text{\rm scal}(p) \leq e^{-\varepsilon \sigma \, d(p_0,p)} \, \text{\rm scal}(p_0)\] 
for $d(p_0,p) \geq 1$. \\

\begin{proposition}[L.~Ma, D.~Chen \cite{Ma-Chen}]
\label{curvature.decay.on.expanding.solitons}
Let $(M,g)$ be an expanding gradient Ricci soliton which is complete and non-compact. Suppose that there exists a point $p_0 \in M$ such that $0 < \text{\rm scal}(p) \leq \text{\rm scal}(p_0)$ for all points $p \in M$. Moreover, we assume that $\text{\rm Ric} \geq \varepsilon \, \text{\rm scal} \: g$ for some constant $\varepsilon > 0$. Then there exists a constant $\alpha > 0$ such that 
\[\text{\rm scal}(p) \leq e^{-\alpha \, d(p_0,p)^2} \, \text{\rm scal}(p_0)\] 
for all points $p \in M$.
\end{proposition}

\textit{Proof of Proposition \ref{curvature.decay.on.expanding.solitons}.} 
Since $(M,g)$ is an expanding gradient Ricci soliton, there exists a smooth function $f: M \to \mathbb{R}$ and a constant $\rho < 0$ such that $\text{\rm Ric}_{ij} = \rho \, g_{ij} + D_{i,j}^2 f$. This implies 
\[\partial_i \text{\rm scal} + 2 \, \text{\rm Ric}_{ij} \, \partial^j f = 0.\] 
By assumption, the scalar curvature attains its maximum at the point $p_0$. This implies $\partial_i \text{\rm scal}(p_0) = 0$. Since $M$ has positive Ricci curvature, it follows that $\partial_i f(p_0) = 0$. Thus, $p_0$ is a critical point of the function $f$. Since $f$ is strictly convex, $p_0$ is the only critical point of $f$. Since $\text{\rm Ric} \geq 0$, we have $f(p) - f(p_0) \geq -\frac{1}{2} \, \rho \, d(p_0,p)^2$ for all $p \in M$. 

As above, the inequality $\text{\rm Ric} \geq \varepsilon \, \text{\rm scal} \: g$ implies 
\[\partial_i(e^{2\varepsilon \, f} \, \text{\rm scal}) \, \partial^i f \leq 0.\] 
Fix an arbitrary point $p \in M$, and let $\gamma: [0,\infty) \to M$ be the solution of the ODE $\gamma'(s) = -\partial^i f(\gamma(s)) \, \frac{\partial}{\partial x^i}$ with initial condition $\gamma(0) = p$. Then $\gamma(s) \to p_0$ as $s \to \infty$. Moreover, the function $s \mapsto e^{2\varepsilon \, f(\gamma(s))} \, \text{\rm scal}(\gamma(s))$ is monotone increasing. Therefore, we have 
\[\text{\rm scal}(p) \leq e^{-2\varepsilon (f(p)-f(p_0))} \, \text{\rm scal}(p_0) \leq e^{\varepsilon \rho \, d(p_0,p)^2} \, \text{\rm scal}(p_0)\] 
for all points $p \in M$. Since $\rho < 0$, the assertion follows. \\

The following theorem generalizes Theorem \ref{chen.zhu.compactness} above: 

\begin{theorem}
\label{generalization.of.chen.zhu}
Let $(M,g_0)$ be a complete Riemannian manifold of dimension $n \geq 4$ with bounded curvature. Suppose that there exists a positive constant $\varepsilon$ such that 
\[R_{1313} + \lambda^2 \, R_{1414} + \mu^2 \, R_{2323} + \lambda^2\mu^2 \, R_{2424} - 2\lambda\mu \, R_{1234} \geq \varepsilon \, \text{\rm scal} > 0\]  
for all orthonormal four-frames $\{e_1,e_2,e_3,e_4\}$ and all $\lambda,\mu \in [-1,1]$. Then $M$ is compact.
\end{theorem}

\textit{Proof of Theorem \ref{generalization.of.chen.zhu}.}
We argue by contradiction. Suppose that $M$ is non-compact. By work of Shi, we can find a maximal solution to the Ricci flow with initial metric $g_0$ (see \cite{Shi}, Theorem 1.1). Let us denote this solution by $g(t)$, $t \in [0,T)$. Using Proposition 13 in \cite{Brendle-Schoen1}, one can show that there exists a positive constant $\delta$ with the following property: for each $t \in [0,T)$, the curvature tensor of $(M,g(t))$ satisfies  
\begin{equation} 
\label{pinching}
R_{1313} + \lambda^2 \, R_{1414} + \mu^2 \, R_{2323} + \lambda^2\mu^2 \, R_{2424} - 2\lambda\mu \, R_{1234} \geq \delta \, \text{\rm scal} 
\end{equation}
for all orthonormal four-frames $\{e_1,e_2,e_3,e_4\}$ and all $\lambda,\mu \in [-1,1]$. The constant $\delta$ depends on $\varepsilon$ and $n$, but not on $t$. In particular, the manifold $(M,g(t))$ has positive sectional curvature for all $t \in [0,T)$. 

By a theorem of Gromoll and Meyer, the injectivity radius of $(M,g(t))$ is bounded from below by 
\[\text{\rm inj}(M,g(t)) \geq \frac{\pi}{\sqrt{N(t)}},\] 
where $N(t) = \sup_{p \in M} \text{\rm scal}_{g(t)}(p)$ denotes the supremum of the scalar curvature of $(M,g(t))$. There are three possibilities: 

\textit{Case 1:} Suppose that $T < \infty$. Let $F$ be a pinching set with the property that the curvature tensor of $g(0)$ lies in $F$ for all points $p \in M$. (The existence of such a pinching set follows from Proposition 17 in \cite{Brendle-Schoen1}.) Using Hamilton's maximum principle for systems, we conclude that the curvature tensor of $g(t)$ lies in $F$ for all points $p \in M$ and all $t \in [0,T)$. 

Since $T < \infty$, we have $\sup_{t \in [0,T)} N(t) = \infty$. Hence, we can find a sequence of times $t_k \in [0,T)$ such that $N(t_k) \to \infty$. Let us dilate the manifolds $(M,g(t_k))$ so that the maximum of the scalar curvature is equal to $1$. These rescaled manifolds converge to a limit manifold $\hat{M}$ which has pointwise constant sectional curvature. Using Schur's lemma, we conclude that $\hat{M}$ has constant sectional curvature. Consequently, $\hat{M}$ is compact by Myers theorem. On the other hand, $\hat{M}$ is non-compact, since it arises as a limit of non-compact manifolds. This is a contradiction. 

\textit{Case 2:} Suppose that $T = \infty$ and $\sup_{t \in [0,\infty)} t \, N(t) = \infty$. By a result of Hamilton, there exists a sequence of dilations of the solution $(M,g(t))$ which converges to a singularity model of Type II (see \cite{Hamilton-survey}, Theorem 16.2). We denote this limit solution by $(\hat{M},\hat{g}(t))$. The solution $(\hat{M},\hat{g}(t))$ is defined for all $t \in (-\infty,\infty)$. Moreover, there exists a point $p_0 \in \hat{M}$ such that 
\[\text{\rm scal}_{\hat{g}(t)}(p) \leq \text{\rm scal}_{\hat{g}(0)}(p_0) = 1\] 
for all points $(p,t) \in \hat{M} \times (-\infty,\infty)$. 

The manifold $(\hat{M},\hat{g}(0))$ satisfies the pinching estimate (\ref{pinching}), as (\ref{pinching}) is scaling invariant. Moreover, it follows from the strict maximum principle that $\text{\rm scal}_{\hat{g}(0)}(p) > 0$ for all $p \in \hat{M}$. Therefore, the manifold $(\hat{M},\hat{g}(0))$ has positive sectional curvature. Since $(\hat{M},\hat{g}(0))$ arises as a limit of complete, non-compact manifolds, we conclude that $(\hat{M},\hat{g}(0))$ is complete and non-compact. By a theorem of Gromoll and Meyer \cite{Gromoll-Meyer}, the manifold $\hat{M}$ is diffeomorphic to $\mathbb{R}^n$.

It follows from Proposition \ref{singularity.model.of.type.2} that $(\hat{M},\hat{g}(0))$ is a steady gradient Ricci soliton. By Proposition \ref{curvature.decay.on.steady.solitons}, the scalar curvature of $(\hat{M},\hat{g}(0))$ decays exponentially. Hence, a theorem of A.~Petrunin and W.~Tuschmann implies that $(\hat{M},\hat{g}(0))$ is isometric to $\mathbb{R}^n$ (see \cite{Petrunin-Tuschmann}, Theorem B). This contradicts the fact that $\text{\rm scal}_{\hat{g}(0)}(p_0) = 1$. 

\textit{Case 3:} Suppose that $T = \infty$ and $\sup_{t \in [0,\infty)} t \, N(t) < \infty$. By a result of Hamilton, there exists a sequence of dilations of the solution $(M,g(t))$ which converges to a singularity model of Type III (see \cite{Hamilton-survey}, Theorem 16.2). We denote this limit solution by $(\hat{M},\hat{g}(t))$. The solution $(\hat{M},\hat{g}(t))$ is defined for all $t \in (-A,\infty)$, where $A$ is a positive real number. Moreover, there exists a point $p_0 \in \hat{M}$ such that 
\[(A + t) \cdot \text{\rm scal}_{\hat{g}(t)}(p) \leq A \cdot \text{\rm scal}_{\hat{g}(0)}(p_0) = A\] 
for all points $(p,t) \in \hat{M} \times (-A,\infty)$. 

As above, the manifold $(\hat{M},\hat{g}(0))$ satisfies the pinching estimate (\ref{pinching}). Moreover, the strict maximum principle implies that $\text{\rm scal}_{\hat{g}(0)}(p) > 0$ for all $p \in \hat{M}$. Consequently, the manifold $(\hat{M},\hat{g}(0))$ has positive sectional curvature. Moreover, the manifold $(\hat{M},\hat{g}(0))$ is complete and non-compact, since it arises as a limit of complete, non-compact manifolds. Therefore, $\hat{M}$ is diffeomorphic to $\mathbb{R}^n$ (see \cite{Gromoll-Meyer}).

By Proposition \ref{singularity.model.of.type.3}, the manifold $(\hat{M},\hat{g}(0))$ is an expanding gradient Ricci soliton. Hence, Proposition \ref{curvature.decay.on.expanding.solitons} implies that the scalar curvature of $(\hat{M},\hat{g}(0))$ decays exponentially. By Theorem B in \cite{Petrunin-Tuschmann}, the manifold $(\hat{M},\hat{g}(0))$ is isometric to $\mathbb{R}^n$. This contradicts the fact that $\text{\rm scal}_{\hat{g}(0)}(p_0) = 1$. \\

This completes the proof of Theorem \ref{generalization.of.chen.zhu}. \\

\begin{corollary} 
Let $(M,g_0)$ be a complete Riemannian manifold of dimension $n \geq 4$ with bounded curvature. Suppose that there exists a positive constant $\varepsilon$ such that $0 < K(\pi_1) < (4-\varepsilon) \, K(\pi_2)$ for all points $p \in M$ and all two-planes $\pi_1,\pi_2 \subset T_p M$. Then $M$ is compact.
\end{corollary}


\begin{thebibliography}{99}
\bibitem{Abresch-Meyer}
U.~Abresch and W.~Meyer, \textit{A sphere theorem with a pinching constant below $1/4$,} J. Diff. Geom. 44, 214--261 (1996)

\bibitem{Andrews-Nguyen}
B.~Andrews and H.~Nguyen, \textit{Four-manifolds with $1/4$-pinched flag curvatures,} to appear in Asian J. Math.

\bibitem{Berger1}
M.~Berger, \textit{Les vari\'et\'es Riemanniennes $1/4$-pinc\'ees,} Ann. Scuola Norm. Sup. Pisa 14, 161--170 (1960)

\bibitem{Berger2}
M.~Berger, \textit{Sur quelques vari\'et\'es riemaniennes suffisamment pinc\'ees,} 
Bull. Soc. Math. France 88, 57--71 (1960)

\bibitem{Berger3}
M.~Berger, \textit{Sur les vari\'et\'es $4/23$-pinc\'ees de dimension $5$,} 
C. R. Acad. Sci. Paris 257, 4122--4125 (1963)

\bibitem{Berger4}
M.~Berger, \textit{Sur les vari\'et\'es riemanniennes pinc\'ees juste au-dessous de $1/4$,} Ann. Inst. Fourier (Grenoble) 33, 135--150 (1983)

\bibitem{Bohm-Wilking} 
C.~B\"ohm and B.~Wilking, \textit{Manifolds with positive curvature operator are space forms,} Ann. of Math. 167, 1079--1097 (2008)

\bibitem{Bony} 
J.M.~Bony, \textit{Principe du maximum, in\'egalit\'e de Harnack et unicit\'e du probl\`eme de Cauchy pour les op\'erateurs elliptiques d\'eg\'en\'er\'es,} Ann. Inst. Fourier (Grenoble) 19, 277--304 (1969) 

\bibitem{Bourguignon}
J.P.~Bourguignon, \textit{La conjecture de Hopf sur $S^2 \times S^2$,} Riemannian geometry in dimension $4$ (Paris 1978/1979), 347--355, Textes Math. 3, CEDIC, Paris (1981) 

\bibitem{Brendle-Schoen1} 
S.~Brendle and R.~Schoen, \textit{Manifolds with $1/4$-pinched curvature are space forms,} J. Amer. Math. Soc. 22, 287--307 (2009)

\bibitem{Brendle-Schoen2}
S.~Brendle and R.~Schoen, \textit{Classification of manifolds with weakly $1/4$-pinched curvatures,} Acta Math. 200, 1--13 (2008)

\bibitem{Brendle1}
S.~Brendle, \textit{A general convergence result for the Ricci flow,} Duke Math. J. 145, 585--601 (2008)

\bibitem{Brendle2} 
S.~Brendle, \textit{A generalization of Hamilton's differential Harnack inequality for the Ricci flow,} J. Diff. Geom. 82, 207--227 (2009)

\bibitem{Brendle3}
S.~Brendle, \textit{Einstein manifolds with nonnegative isotropic curvature are locally symmetric,} to appear in Duke Math. J.

\bibitem{Cao} 
H.D.~Cao, \textit{On Harnack's inequalities for the K\"ahler-Ricci flow,} Invent. Math. 109, 247--263 (1992)

\bibitem{Chang-Gursky-Yang}
A.~Chang, M.~Gursky, and P.~Yang, \textit{A conformally invariant sphere theorem in four dimensions,} Publ. Math. IH\'ES 98, 105--143 (2003)

\bibitem{Cheeger-Gromoll}
J.~Cheeger and D.~Gromoll, \textit{The splitting theorem for 
manifolds of nonnegative Ricci curvature,} J. Diff. Geom. 6, 119--128 (1971)

\bibitem{Chen-Zhu1}
B.~Chen and X.~Zhu, \textit{Complete Riemannian manifolds with pointwise pinched curvature,} Invent. Math. 423--452 (2000)

\bibitem{Chen-Zhu2} 
B.~Chen and X.~Zhu, \textit{Ricci flow with surgery on four-manifolds with positive isotropic curvature,} J. Diff. Geom. 74, 177--264 (2006)

\bibitem{Chen} 
H.~Chen, \textit{Pointwise $1/4$-pinched $4$-manifolds,} Ann. Global Anal. Geom. 9, 161--176 (1991)

\bibitem{DeTurck}
D.~DeTurck, \textit{Deforming metrics in the direction of their Ricci tensors,} J. Diff. Geom. 18, 157--162 (1983)

\bibitem{Fraser1}
A.~Fraser, \textit{On the free boundary variational problem for minimal disks,} Comm. Pure Appl. Math. 53, 931--971 (2000)

\bibitem{Fraser2}
A.~Fraser, \textit{Fundamental groups of manifolds with positive isotropic curvature,} Ann. of Math. 158, 345--354 (2003)

\bibitem{Greene-Wu1}
R.E.~Greene and H.~Wu, \textit{On the subharmonicity and plurisubharmonicity of geodesically 
convex functions,} Indiana Univ. Math. J. 22, 641--653 (1972)

\bibitem{Greene-Wu2}
R.E.~Greene and H.~Wu, \textit{$C^\infty$ convex functions and manifolds of positive curvature,} 
Acta Math. 137, 209--245 (1976)

\bibitem{Gromoll} 
D.~Gromoll, \textit{Differenzierbare Strukturen und Metriken positiver Kr\"ummung auf Sph\"aren,} Math. Ann. 164, 353--371 (1966)

\bibitem{Gromoll-Grove}
D.~Gromoll and K.~Grove, \textit{A generalization of Berger's rigidity theorem for positively curved manifolds,}
Ann. Sci. \'Ecole Norm. Sup. 20, 227--239 (1987)

\bibitem{Gromoll-Meyer} 
D.~Gromoll and W.~Meyer, \textit{On complete open manifolds of positive curvature,} Ann. of Math. 20, 75--90 (1969)

\bibitem{Gromov}
M.~Gromov, \textit{Positive curvature, macroscopic dimension, spectral gaps and higher signatures,} Functional analysis on the eve of the 21st century, Vol. II (New Brunswick 1993), 1--213,
Progr. Math., 132, Birkh\"auser, Boston (1996)

\bibitem{Grothendieck}
A.~Grothendieck, \textit{Sur la classification des fibr\'es holomorphes sur la sph\`ere de Riemann,} Amer. J. Math. 79, 121--138 (1957)

\bibitem{Grove-Karcher-Ruh1}
K.~Grove, H.~Karcher, and E.~Ruh, \textit{Group actions and curvature,}  Invent. Math. 23, 
31--48 (1974)

\bibitem{Grove-Karcher-Ruh2}
K.~Grove, H.~Karcher, and E.~Ruh, \textit{Jacobi fields and Finsler metrics on compact Lie groups with an application to differentiable pinching problems,} Math. Ann. 211, 7--21 (1974)

\bibitem{Grove-Shiohama}
K.~Grove and K.~Shiohama, \textit{A generalized sphere theorem,} Ann. of Math. 106, 201--211 (1977)

\bibitem{Grove-survey}
K.~Grove, \textit{Ramifications of the classical sphere theorem,} Actes de la Table Ronde de G\'eom\'etrie Diff\'erentielle (Luminy 1992), 363--376, S\'emin. Congr. 1, Soc. Math. France (1996)

\bibitem{Gulliver}
R.~Gulliver, \textit{Regularity of minimizing surfaces of prescribed mean curvature,} Ann. of Math. 97, 275--305 (1973)

\bibitem{Hamilton1}
R.~Hamilton, \textit{Three-manifolds with positive Ricci curvature,} J. Diff. Geom. 17, 255--306 (1982)

\bibitem{Hamilton2} 
R.~Hamilton, \textit{Four-manifolds with positive curvature operator,} J. Diff. Geom. 24, 153--179 (1986)

\bibitem{Hamilton3}
R.~Hamilton, \textit{The Ricci flow on surfaces,} Contemp. Math. 71, 237--262 (1988)

\bibitem{Hamilton4} 
R.~Hamilton, \textit{The Harnack estimate for the Ricci flow,} J. Diff. Geom. 37, 225--243 (1993)

\bibitem{Hamilton5} 
R.~Hamilton, \textit{Eternal solutions to the Ricci flow,} J. Diff. Geom. 38, 1--11 (1993)

\bibitem{Hamilton6}
R.~Hamilton, \textit{Convex hypersurfaces with pinched second fundamental form,} Comm. Anal. Geom. 2, 167--172 (1994)

\bibitem{Hamilton7} 
R.~Hamilton, \textit{A compactness property for solutions of the Ricci flow,} Amer. J. Math. 117, 545--572 (1995)

\bibitem{Hamilton-survey}
R.~Hamilton, \textit{The formation of singularities in the Ricci flow,} Surveys in Differential Geometry 2, 7--136 (1995)

\bibitem{Hamilton8} 
R.~Hamilton, \textit{Four-manifolds with positive isotropic curvature,} Comm. Anal. Geom. 5, 1--92 (1997)

\bibitem{Huisken} 
G.~Huisken, \textit{Ricci deformation of the metric on a Riemannian manifold,} J. Diff. Geom. 21, 47--62 (1985)

\bibitem{Im-Hof-Ruh}
H.~Im~Hof and E.~Ruh, \textit{An equivariant pinching theorem,} Comment. Math. Helv. 50, no. 3, 
389--401(1975)

\bibitem{Karcher} 
H.~Karcher, \textit{A short proof of Berger's curvature tensor estimates,} Proc. Amer. Math. Soc. 26, 642--644 (1970)

\bibitem{Klingenberg}
W.~Klingenberg, \textit{\"Uber Riemannsche Mannigfaltigkeiten mit positiver Kr\"ummung,} Comment. Math. Helv. 35, 47--54 (1961)

\bibitem{Ma-Chen}
L.~Ma and D.~Chen, \textit{Remarks on non-compact complete Ricci expanding solitons,} arxiv:0508363

\bibitem{Margerin1}
C.~Margerin, \textit{Pointwise pinched manifolds are space forms,} Geometric measure theory and the calculus of variations (Arcata 1984), 343--352, Proc. Sympos. Pure Math. 44, Amer. Math. Soc., Providence RI (1986)

\bibitem{Margerin2}
C.~Margerin, \textit{A sharp characterization of the smooth $4$-sphere in curvature terms,} Comm. Anal. Geom. 6, 21--65 (1998)

\bibitem{Micallef-Moore}
M.~Micallef and J.D.~Moore, \textit{Minimal two-spheres and the topology of manifolds with positive curvature on totally isotropic two-planes,} Ann. of Math. 127, 199-227 (1988)

\bibitem{Micallef-Wang}
M.~Micallef and M.~Wang, \textit{Metrics with nonnegative isotropic curvature,} Duke Math. J. 72, no. 3, 649--672 (1993)

\bibitem{Micallef-White}
M.~Micallef and B.~White, \textit{The structure of branch points in minimal surfaces and in pseudoholomorphic curves,} Ann. of Math. 141, 35--85 (1995)

\bibitem{Nguyen}
H.~Nguyen, \textit{Invariant curvature cones and the Ricci flow,} PhD thesis, Australian National University (2007)

\bibitem{Ni}
L.~Ni, \textit{Ancient solutions to K\"ahler-Ricci flow,} Math. Res. Lett. 12, 633--653 (2005)

\bibitem{Ni-Wu}
L.~Ni and B.~Wu, \textit{Complete manifolds with nonnegative curvature operator,} Proc. Amer. Math. Soc. 135, 3021--3028 (2007)

\bibitem{Nishikawa}
S.~Nishikawa, \textit{Deformation of Riemannian metrics and manifolds with bounded curvature ratios,} Geometric measure theory and the calculus of variations (Arcata 1984), 343--352, Proc. Sympos. Pure Math. 44, Amer. Math. Soc., Providence RI (1986)

\bibitem{Perelman1} 
G.~Perelman, \textit{The entropy formula for the Ricci flow and its geometric applications,} arxiv:0211159

\bibitem{Perelman2}
G.~Perelman, \textit{Ricci flow with surgery on three-manifolds,} arxiv:0303109

\bibitem{Petersen}
P.~Petersen, \textit{Riemannian Geometry,} Graduate Texts in Mathematics, vol. 171, 2nd edition, Springer-Verlag, New York, 2006

\bibitem{Petersen-Tao}
P.~Petersen and T.~Tao, \textit{Classification of almost quarter-pinched manifolds,} Proc. Amer. Math. Soc. 137, 2437--2440 (2009)

\bibitem{Petrunin-Tuschmann}
A.~Petrunin and W.~Tuschmann, \textit{Asymptotical flatness and cone structure at infinity,} Math. Ann. 321, 775--788 (2001)

\bibitem{Rauch} 
H.E.~Rauch, \textit{A contribution to differential 
geometry in the large,} Ann. of Math. 54, 38--55 (1951)

\bibitem{Ruh1} 
E.~Ruh, \textit{Kr\"ummung und differenzierbare Struktur auf Sph\"aren II,} Math. Ann. 205, 113--129 (1973) 

\bibitem{Ruh2}
E.~Ruh, \textit{Riemannian manifolds with bounded curvature ratios,} J. Diff. Geom. 17, 643--653 (1982)

\bibitem{Sacks-Uhlenbeck}
J.~Sacks and K.~Uhlenbeck, \textit{The existence of minimal immersions of $2$-spheres,} Ann. of Math. 113, 1--24 (1981)

\bibitem{Schoen-Yau}
R.~Schoen and S.T.~Yau, \textit{Existence of incompressible minimal surfaces and the topology of three dimensional
manifolds with non-negative scalar curvature,} Ann. of Math. 110, 127--142 (1979)

\bibitem{Schoen-survey}
R.~Schoen, \textit{Minimal submanifolds in higher codimensions,} Mat. Contemp. 30, 169--199 (2006)

\bibitem{Seaman}
W.~Seaman, \textit{A pinching theorem for four manifolds,} Geom. Dedicata 31, 37--40 (1989)

\bibitem{Shi}
W.X.~Shi, \textit{Deforming the metric on complete Riemannian manifolds,} J. Diff. Geom. 30, 223--301 (1989)

\bibitem{Sugimoto-Shiohama} 
M.~Sugimoto and K.~Shiohama, and H.~Karcher, \textit{On the differentiable pinching problem,} Math. Ann. 195, 1--16 (1971)

\bibitem{Tachibana}
S.~Tachibana, \textit{A theorem on Riemannian manifolds with positive curvature operator,} Proc. Japan Acad. 50, 301--302 (1974)

\bibitem{Wilking}
B.~Wilking, \textit{Index parity of closed geodesics and rigidity of Hopf fibrations,} Invent. Math. 144, 281--295 (2001)
\end{thebibliography}
\end{document}